\documentclass[10pt]{amsart}

\usepackage[margin=1.5in]{geometry}
\usepackage{colortbl}
\newtheorem{theorem}{Theorem}[section]

\theoremstyle{definition}

\theoremstyle{remark}
\newtheorem{remark}[theorem]{Remark}

\numberwithin{equation}{section}

\usepackage{graphicx}
\usepackage{subcaption}
\usepackage{multirow}
\usepackage{bbm}
\usepackage{amssymb}

\begin{document}

\title[Conservative Semi-Lagrangian schemes for a BGK model of inert gas mixtures]{Conservative Semi-Lagrangian schemes for a general consistent BGK model for inert gas mixtures}

\author[S.Y. Cho]{Seung Yeon Cho}
\address{Seung Yeon Cho\\
	Department of Mathematics and Computer Science\\
	University of Catania\\
	95125 Catania, Italy}
\email{chosy89@skku.edu}

\author[S. Boscarino]{Sebastiano Boscarino}
\address{Sebastiano Boscarino\\
Department of Mathematics and Computer Science\\
University of Catania\\
95125 Catania, Italy} \email{boscarino@dmi.unict.it}

\author[M. Groppi]{Maria Groppi}
\address{Maria Groppi\\
Department of Mathematical, Physical and Computer Sciences, University of Parma\\ Parco Area delle Scienze 53/A,
43124 Parma, Italy} \email{maria.groppi@unipr.it}

\author[G. Russo]{Giovanni Russo}
\address{Giovanni Russo\\
	Department of Mathematics and Computer Science\\
	University of Catania\\
	95125 Catania, Italy} \email{russo@dmi.unict.it}


\begin{abstract}
In this work, we propose a class of high order semi-Lagrangian scheme for a general consistent BGK model for inert gas mixtures. The proposed scheme not only fulfills indifferentiability principle, but also asymptotic preserving property, which allows us to capture the behaviors of hydrodynamic limit models. We consider two hydrodynamic closure which can be derived from the BGK model at leading order: classical Euler equations for number densities, global velocity and temperature, and a multi-velocities and temperatures Euler system. Numerical simulations are performed to demonstrate indifferentiability principle and asymptotic preserving property of the proposed conservative semi-Lagrangian scheme to the Euler limits.
\end{abstract}
\maketitle
\section{Introduction}

Relaxation time approximations of BGK type for the Boltzmann equation are fundamental tools in order to proper describe several  regimes
of gas dynamics, where the macroscopic approach is not adequate. Indeed, the Boltzmann equation \cite{Cer} is difficult to treat both from an analytical and a numerical point of view, due to its highly nonlinear integral collision operator, whereas BGK models \cite{BGK, W} retain the essential features of the Boltzmann equation with a simpler structure. The extension of BGK models to mixtures of gases is not trivial, since some drawbacks can arise, like loss of positivity of temperature and violation of the indifferentiability principle.
To our knowledge, the first  mathematical approach to the construction of BGK models for mixtures was presented  in \cite{AAP}, where the questions of consistency of the model and its generality were addressed. Generality means that the model is valid for any number of components of the mixture; consistency means exact conservation laws, fulfilment of the H-Theorem and uniqueness of equilibrium solutions, The model proposed in \cite{AAP} was based on a single relaxation operator for each species, which reproduce the structure of the BGK operator for a single gas. Later, different BGK models have been introduced, having the same structure and allowing also to describe simple bimolecular chemical reactions in mixtures \cite{PoFGS, GRS}  as well as polyatomic gases \cite{BC}.

Recently, in \cite{BBGSP1} a general consistent BGK model for inert mixtures of monoatomic gases has been proposed, which mimics instead the structure of the Boltzmann equations for mixtures, with a collision operator for each species that is a sum of bi-species BGK operators.

In this paper we propose high order semi-	Lagrangian schemes for this general BGK model for mixtures. Semi-Lagrangian schemes have been widely adopted for solving BGK type-models of the Boltzmann equation, since they do not require stability restriction on time step. In \cite{RS}, authors first applied semi-Lagrangian method to the classical BGK equation for a single gas, combined  with L-stable time discretization, which enables one to avoid the time step restrictions coming from the convective term and small relaxation time. This approach has been then adopted with linear multi-step methods such as backward difference formula (BDF) in \cite{GRS} to design more efficient schemes. 	
In \cite{BCRY}, authors verify that the schemes \cite{GRS,RS} are not conservative and cannot exactly capture the behavior of Euler system when the Knudsen number $\varepsilon$ tends to zero. This so called asymptotic preserving property (AP) can be attained indeed using conservative schemes, (see also \cite{JX}). To this aim, a predictor-corrector type semi-Lagrangian scheme is introduced in \cite{BCRY} to attain conservation up to machine precision. Although this AP scheme allows one to use large time step, it has severe time step restriction for large relaxation time. As a remedy of this, in \cite{BCRY3} authors also modify the high order semi-Lagrangian methods \cite{RS,GRS} by making use of a conservative reconstruction technique \cite{BCRY2} for the interpolation together with the construction of conservative Maxwellian introduced in \cite{M,GT}. This conservative scheme allows large time step for all ranges of relaxation time, and capture the correct behavior of Euler system in the limit $\varepsilon \to 0$. 
In the context of inert gas mixtures, a semi-Lagrangian approach has been already applied to BGK-type models in \cite{GRS2}. Although that method is designed to be high order without splitting based on \cite{GRS}, it is not conservative. 
In this paper, we focus on high order conservative semi-Lagrangian schemes with asymptotic property and we apply them to a general consistent BGK model for inert mixture \cite{BBGSP1}, making use also of conservative reconstruction technique in \cite{BCRY2}.

The semi-Lagrangian approximation of the BGK equations is tested versus the Euler systems that can be obtained as a zero order level of any asymptotic expansion with respect to a proper Knudsen number. In this respect two hydrodynamic closures can be deduced from general consistent model \cite{BBGSP1}, leading to classic Euler equations for number densities, global velocity and temperature or to multi-velocity and multi-temperature Euler system. 

The paper is organized as follows. In section 2, we review the basic properties of the Boltzmann equations, the general BGK model for mixtures and their macroscopic Euler closures. Section 3 is devoted to the derivation of the first order semi-Lagrangian scheme for the general BGK model. In section 4, we extend it to high order semi-Lagrangian schemes. Then, in section 5 we present several numerical tests to verify the performance of the proposed schemes, with particular attention to the indifferentiability principle and AP property. Some concluding remarks end the paper.

\section{Kinetic models for inert gas mixtures}

\subsection{The Boltzmann equation for inert gas mixtures.}
Let us consider distribution functions $\left\{ f_s(x,\text{v},t),~ s=1,\cdots,L \right\}$ of $L$-species inert gases defined on phase space $(x,v)\in \mathbb{R}^3 \times \mathbb{R}^3$ at time $t>0$. Under the assumption that $s$-th gas has a mass $m_s>0$,
the dynamics of a mixture of $L$ inert rarefied gases can be described by the set of $L$ Boltzmann-type equations:
\begin{align}\label{Boltzmann eqn}
	\frac{\partial {f_s}}{\partial t} + \text{v}\cdot \nabla_x f_s = Q_{s},\quad s=1,\cdots,L,
\end{align}
where the collision operator $Q_s$ for $s$-species gas is given by 
\[
Q_s=\sum_{k=1}^L Q_{sk}(f_s,f_k).
\]
The binary collision operator $Q_{sk}$ is defined by
\begin{align}\label{Qsk}
	Q_{sk}(f_s,f_k)= \int_{\mathbb{R}^3 \times \mathbb{S}^2} d\text{w}\,d\omega\, g_{sk}(|y|,\hat{y}\cdot \omega) \bigg[ f_s(\text{v}')f_k(\text{w}')-f_s(\text{v})f_k(\text{w})\bigg],
\end{align}
where $g_{sk}$ is a nonnegative scattering kernel. In \eqref{Qsk}, we use an integration variable $\text{w} \in \mathbb{R}^3$, a unit vector on a sphere $\omega \in \mathcal{S}^2$, a relative velocity $y:=\text{v}-\text{w}$ and its unit vector $\hat{y}:=y/|y|$. In the collision between $s$- and $k$-species gases, the pre-collisional velocities $\text{v}$, $\text{w}$ and post collisional velocities $\text{v}'$ and $\text{w}'$ are related as follows:
\[
\text{v}'= \frac{m_s \text{v} + m_k\text{w}}{m_s + m_k} + \frac{m_k}{m_s + m_k}|y|\omega,
\quad 
\text{w}'= \frac{m_s \text{v} + m_k\text{w}}{m_s + m_k} - \frac{m_s}{m_s + m_k}|y|\omega.
\]
We recall that the collision operator $Q_s$ has the following properties:
\begin{itemize}
	\item Conservation laws:
	\[
	\langle Q_s, 1\rangle=0, \quad \sum_{s=1}^L m_s\langle Q_s, \text{v}\rangle=0, \quad \sum_{s=1}^L m_s\langle Q_s, |\text{v}|^2\rangle=0
	\]
	where
	\begin{align*}
	\langle f, h\rangle := \int_{\mathbb{R}^3} dv f(\text{v})h(\text{v}).
	\end{align*}
	\item H-theorem:
	\[
	\sum_{s=1}^L \langle Q_s , \log f_s\rangle \leq 0
	\]
	\item Uniqueness of equilibrium:
	\[Q_s=0, \quad s=1,\cdots, L \qquad \Longrightarrow  \qquad f_s= n_s M(\text{v};u_s,\frac{K_BT_s}{m_s})
	\]
	where 
	\begin{align}\label{M form}
		M=M(\text{v};a,b) \equiv \left(\frac{1}{2\pi b}\right)^{3/2} \exp \left(- \frac{1}{2b}\left|\text{v}-a\right|^2\right)
	\end{align}
	and $K_B$ stands for the Boltzmann constant.
	\item Indifferentiability principle: this means that the sum of distribution function $f=\sum_{s=1}^L f_s$ obeys the single gas Boltzmann equation if all $m_s$ and collision kernel $g_{sk}$ are identical.

	\item Exchange rates of momentum and energy for each pair of species $s$ and $k$:
	\begin{align*}
	\langle Q_{sk}, \text{v} \rangle &= -\frac{m_k}{m_s + m_k} \langle\langle \text{v}- \text{w} \rangle\rangle_{sk}\cr
	\langle Q_{sk}, |\text{v}|^2 \rangle &= -\frac{2m_k}{(m_s + m_k)^2} \langle\langle (m_s\text{v} + m_k\text{w}) \cdot (\text{v}- \text{w}) \rangle\rangle_{sk}
	\end{align*}
	where 
	\begin{align*}
	\langle\langle \psi(\text{v},\text{w}) \rangle\rangle_{sk} &= \int_{\mathbb{R}^3 \times \mathbb{R}^3} d\text{v}\, d\text{w} f_s(\text{v}) f_k(\text{w}) g_{sk}^{(1)}(|\text{v}-\text{w}|) \psi(\text{v},\text{w})\cr
	g_{sk}^{(1)}(|y|)&= 2 \pi \int_{-1}^{1} d\mu (1-\mu) g_{sk}(|y|,\mu), \quad y \in \mathbb{R}^3.
	\end{align*}	
	In this paper, we consider Maxwell molecules, for which the cross section $g_{sk}$ is independent of $|y|$, i.e., \begin{align}\label{lambda derivation}
		g_{sk}^{(1)}(|y|)=\lambda_{sk}=\text{const},
	\end{align}
	Consequently,
	\begin{align}\label{ex rate}
	\begin{split}
		\langle Q_{sk}, \text{v} \rangle &= -\frac{m_k}{m_s + m_k} \lambda_{sk}n_sn_k(u_s-u_k)\cr
		\langle Q_{sk}, |\text{v}|^2 \rangle &= -\frac{2m_k}{(m_s + m_k)^2} \lambda_{sk}n_sn_k \left[ 3(K_BT_s-K_BT_k) + (m_su_s + m_ku_k) \cdot (u_s-u_k)\right]
	\end{split}
	\end{align}
	For the details of this calculation and other possible choices of $g_{sk}$, we refer to \cite{BBGSP1} and references therein.
	
	\item Macroscopic variables:
	for $s$-species gas, we define number density $n_s$, average velocity $u_s$, absolute temperature $T_s$ by
	\begin{align}\label{def n nu T}
	n_s=\langle f_s,1 \rangle, \quad n_s u_s=\langle f_s,\text{v} \rangle, \quad 3n_sK_BT_s=m_s\langle f_s,|\text{v}-u_s|^2 \rangle 
	\end{align}
    In a similar way, global macroscopic variables such as number density $n$, density $\rho$, velocity $u$ and temperature $T$ of the whole mixture are given as follows:
	\begin{align*}
	n=\sum_{s=1}^L n_s,\quad \rho=\sum_{s=1}^L \rho_s,\quad \rho_s = m_s n_s,\quad s=1,\cdots,L\cr
	u= \frac{1}{\rho} \sum_{s=1}^L \rho_s u_s, \quad 3nK_BT=3\sum_{s=1}^L n_sK_BT_s + \sum_{s=1}^L \rho_s|u_s-u|^2	
	\end{align*}
\end{itemize}

\subsection{A general consistent BGK model for inert gas mixtures}
In \cite{BBGSP1}, authors introduce a BGK-type model (BBGSP model) which  prescribes the same exchange rates \eqref{ex rate} of the binary collision operators $Q_{sk}$ of Boltzmann equation \eqref{Boltzmann eqn}. In detail, this model uses a relaxation towards fictitious Maxwellian $M_{sk}$ in the place of the binary collision operator $Q_{sk}$ of Boltzmann equation. The model equations are given by 
\begin{align}\label{bgk bbgsp}
\frac{\partial{f_s}}{\partial{t}} + \text{v} \cdot \nabla_x{f_s} =  \frac{1}{\varepsilon}\sum_{k=1}^{L} \nu_{sk}\left(n_{s}M_{sk}-f_s\right),\quad s=1,\cdots, L
\end{align}
where $\varepsilon$ is the Knudsen number, $\nu_{sk}>0$ are free parameters and  $M_{sk}:=M(\text{v};u_{sk},\frac{K_BT_{sk}}{m_s})$ are auxiliary Maxwellians depending on fictitious parameter $u_{sk}$, $T_{sk}$, which are not the true moments of the distribution functions. In order to guarantee the same exchange rates of the Boltzmann equations, the parameters
$u_{sk}$ and $T_{sk}$ take the following forms:
\begin{align}
\begin{split}\label{uskTsk}
	u_{sk}&= (1-a_{sk})u_s + a_{sk}u_k\cr
	T_{sk}&= (1-b_{sk})T_s + b_{sk}T_k +  \frac{\gamma_{sk}}{K_B}\left|u_{s}-u_{k}\right|^2
\end{split}
\end{align}
with
\begin{align}\label{abgamma original}
	a_{sk}= \frac{\lambda_{sk} n_k m_k}{\nu_{sk}(m_s + m_k)}, \quad 
	b_{sk}= \frac{2a_{sk} m_s}{m_s + m_k}, \quad \gamma_{sk}= \frac{m_s a_{sk}}{3}\left( \frac{2m_k}{m_s + m_k} -a_{sk}\right)
\end{align}
In spite of the simpler approximate relaxation operator, the model \eqref{bgk bbgsp} still satisfies qualitative properties of Boltzmann equation such as Conservation law, $H$-theorem, indifferentiability principle. One thing we should keep in mind is that the positivity of $T_{sk}$ should be guaranteed for this model to be well-defined. For this, as remarked in \cite{BBGSP1}, the collision frequency $\nu_{sk}$ should satisfy
\begin{align}\label{positivity condition}
	\nu_{sk} \geq \frac{1}{2}\lambda_{sk} n_k, \quad \text{given}\quad T_s,T_k>0.
\end{align}


\subsubsection{Euler limit} 
In \cite{BBGSP2}, the Chapman Enskog asymptotic expansion has been formally applied to BGK equations \eqref{bgk bbgsp} in order to obtain proper hydrodynamic closures of the macroscopic equations. In particular, at zero-th order approximation in the collision dominated regime the following classical Euler system can be obtained: 
\begin{align}\label{Euler}
\begin{split}
&\frac{\partial n_s}{\partial t} + \nabla_x \cdot(n_s u) =0, \quad s=1,\cdots,L\cr
&\frac{\partial}{\partial t}(\rho u) + \nabla_x \cdot(\rho u \otimes u) + \nabla_x(nK_B T)=0,\cr
&\frac{\partial}{\partial t}\left(\frac{1}{2} \rho|u|^2 + \frac{3}{2}nK_BT\right) + \nabla_x \cdot\left[\left(\frac{1}{2}\rho |u|^2 + \frac{5}{2}nK_BT\right)u\right]=0.
\end{split}
\end{align}

\subsubsection{Multi-temperature and velocity Euler limit.}
The structure of the BGK model allows to consider also a hydrodynamic limit where only intra-species collisions are dominant while other inter-species collisions occurs less often. This approach has been described in \cite{BGM}. We consider the following rescaled form of \eqref{bgk bbgsp}:	
\begin{align}\label{bgk bbgsp multi}
\frac{\partial{f_s}}{\partial{t}} + \text{v} \cdot \nabla_x{f_s} =  \frac{1}{\varepsilon} \nu_{ss}\left(n_{s}M_{ss}-f_s\right) + \frac{1}{\kappa} \sum_{k\neq s}^{L}  \nu_{sk}\left(n_{s}M_{sk}-f_s\right).
\end{align}
where $\kappa$ is a positive parameter; it is  possible to derive hydrodynamic equations for the macroscopic moments of each species $n_s, u_s, T_s$.
Here we study the behavior of Euler system obtained from \eqref{bgk bbgsp multi} at zero order when $ \epsilon \to 0$:
\begin{align}\label{Euler multi}
\begin{split}
&\frac{\partial n_s}{\partial t} + \nabla_x \cdot(n_s u_s) =  0,\cr
&\frac{\partial}{\partial t}(\rho_s u_s) + \nabla_x \cdot(\rho_s u_s \otimes u_s) + \nabla_x(n_sK_BT_s) = \frac{1}{\kappa}\sum_{k\neq s }^{L}\mathcal{R}_{sk},\cr
&\frac{\partial}{\partial t}\left(\frac{1}{2} \rho_s|u_s|^2 + \frac{3}{2}n_sK_BT_s\right) + \nabla_x \cdot\left[\left(\frac{1}{2}\rho_s |u_s|^2 + \frac{5}{2}n_sK_BT_s\right)u_s\right]=\frac{1}{\kappa}\sum_{k\neq s }^{L}\mathcal{S}_{sk}. 
\end{split}
\end{align}
for $s=1,\cdots,L$. The relaxation terms $\mathcal{R}_{sk}$ and $\mathcal{S}_{sk}$ on the r.h.s. are given by
\begin{align}\label{RS form}
\begin{split}
\mathcal{R}_{sk}&= \lambda_{sk}m_{sk}n_sn_k(u_k-u_s),\cr 
\mathcal{S}_{sk}&= \lambda_{sk}\frac{m_{sk}}{m_s + m_k}n_sn_k\left[(m_su_s+m_ku_k)\cdot (u_k-u_s) + 3K_B(T_k - T_s)\right],
\end{split}
\end{align}
where $\displaystyle m_{sk} = \frac{m_sm_k}{m_s + m_k}$.
Note that each gas may have different velocities and temperatures. 
However, in the limit $\kappa \rightarrow 0$, we can expect that the multi-velocities and temperatures solution of \eqref{Euler multi} will converge to the single velocity and temperature solution of \eqref{Euler}. This behavior will be checked in Section \ref{sec numerical test}.
\begin{remark}
	Multi-temperature Euler equations can be obtained also in the framework of Extended Thermodynamics \cite{Simic1,Simic2}. Furthermore, the multi-temperature description is able to better reproduce experimental results for mixtures of noble gases, like for instance the shock structure in Helium-Argon mixtures \cite{Simic3}.
\end{remark}

\subsection{Chu reduction}
For the solution of the kinetic BGK equations \eqref{bgk bbgsp multi}, it is possible to treat problems in slab geometry, as suitable problems which are 1D in both velocity and space. Indeed, the idea of Chu reduction in \cite{chu} enables one to reduce each of the (6+1)-dimensional equation of \eqref{bgk bbgsp multi} to a pair of two (2+1)-dimensional equations. The system can be derived in case of domains with an
axial symmetry with respect to an axis (say $x$-axis). In such a case, along the other two directions ($y$ and $z$ directions) all gradient values of distribution functions and macroscopic velocities $u_2,\,u_3$ vanish. Considering this and taking moments of \eqref{bgk bbgsp multi} with respect to $(1, \text{v}_2^2+\text{v}_3^3)$, one can obtain the following system:
\begin{align}\label{chu bgk bbgsp}
	\begin{split}
		\partial_t g_{1}^{s} + v \partial_x g_1^s&= \frac{1}{\varepsilon}  \nu_{ss} \left(n_{s}M_{ss,1}-g_{1}^{s}\right) + \frac{1}{\kappa} \sum_{k\neq s}^{L} \nu_{sk} \left(n_{s}M_{sk,1}-g_{1}^{s}\right),\quad M_{sk,1}=M(v;u_{sk},\frac{K_BT_{sk}}{m_s})\\
		\partial_t g_{2}^{s} + v \partial_x g_2^s&= \frac{1}{\varepsilon} \nu_{ss} \left(n_{s}M_{ss,2}-g_{2}^{s}\right) + \frac{1}{\kappa} \sum_{k\neq s}^{L} \nu_{sk} \left(n_{s}M_{sk,2}-g_{2}^{s}\right),\quad M_{sk,2}=\frac{2K_BT_{sk}}{m_s}M_{sk,1},
	\end{split}
\end{align}
where
\begin{align*}
	g_1^s(t,x,v)= \int_{\mathbb{R}} f(t,x,\text{v}) d\text{v}_2d\text{v}_3, \quad g_2^s(t,x,v)= \int_{\mathbb{R}} (\text{v}_2^2+\text{v}_3^2)f(t,x,\text{v}) d\text{v}_2d\text{v}_3
\end{align*}
Here the first component of velocity variable is denoted by  $v\equiv\text{v}_1$. We remark that the macroscopic quantities in \eqref{chu bgk bbgsp} can be computed as in \eqref{def n nu T}. Using
\begin{align*}
	n_s= \int_{\mathbb{R}} g_1^s dv,\quad n_s u_s=\int_{\mathbb{R}} v g_1^s dv,\quad \frac{3n_sK_BT_s}{2m_s}=\int_{\mathbb{R}} \frac{|v-u_s|^2}{2}g_1^s  + \frac{g_2^s}{2}dv,
\end{align*}
with \eqref{uskTsk} and \eqref{abgamma original}, one can compute $u_{sk}$ and $T_{sk}$ from $g_1^s$ and $g_2^s$.

	Based on the Chu reduction, our numerical schemes and simulations will be focused on the one-dimensional problem in space. To avoid confusion, we end this section with the representation of one-dimensional hydrodynamic models we will use in the rest of this paper:
\begin{enumerate}
	\item One-dimensional single velocity and temperature Euler system for \eqref{Euler}:
	\begin{align}\label{Euler 1d}
		\begin{split}
			&\frac{\partial n_s}{\partial t} + \partial_x (n_s u) =0, \quad s=1,\cdots,L\cr
			&\frac{\partial}{\partial t}(\rho u) + \partial_x (\rho u^2+n K_BT)=0,\cr
			&\frac{\partial}{\partial t}\left(\frac{1}{2} \rho|u|^2 + \frac{3}{2}nK_BT\right) + \partial_x \left[\left(\frac{1}{2}\rho |u|^2 + \frac{5}{2}nK_BT\right)u\right]=0.
		\end{split}
	\end{align}
	\item One-dimensional multi velocities and temperatures Euler system for \eqref{Euler multi}:
	\begin{align}\label{Euler multi 1d}
		\begin{split}
			&\frac{\partial n_s}{\partial t} + \partial_x (n_s u_s) =  0,\cr
			&\frac{\partial}{\partial t}(\rho_s u_s) + \partial_x (\rho_s u_s^2 + n_s K_BT_s) = \frac{1}{\kappa}\sum_{k\neq s }^{L}\mathcal{R}_{sk},\cr
			&\frac{\partial}{\partial t}\left(\frac{1}{2} \rho_s|u_s|^2 + \frac{3}{2}n_sK_BT_s\right) + \partial_x \left[\left(\frac{1}{2}\rho_s |u_s|^2 + \frac{5}{2}n_sK_BT_s\right)u_s\right]=\frac{1}{\kappa}\sum_{k\neq s }^{L}\mathcal{S}_{sk},
		\end{split}
	\end{align}
	for $s=1,\cdots,L$. Here $\mathcal{R}_{sk}$ and $\mathcal{S}_{sk}$ are computed as in \eqref{RS form}
\end{enumerate}


\section{Semi-Lagrangian methods for a general consistent BGK model}
In this section, we propose a semi-Lagrangian method (SL) in the finite difference framework for solving \eqref{bgk bbgsp multi}. The SL approach has been adopted for various BGK-type models in the context of single gas BGK model \cite{BCRY,BCRY3,GRS,RS} or ES-BGK type models \cite{BCRY4,RY} and a gas mixture model \cite{GRS2}. The main advantage of the technique is to avoid both CFL-type time step restriction of convection term and the stiffness arising in the hydrodynamic limit $\varepsilon \rightarrow 0$.

In the following description, we assume a fixed time step $\Delta t$ and uniform mesh sizes $\Delta x$, $\Delta v$ for space and velocity, respectively. The corresponding grid points will be denoted by 
\begin{align*}
	t_n &= n \Delta t,\qquad\qquad\qquad\quad n=0,1,\cdots,N_t\cr
	v_j &= v_{min} + (j-1)\Delta v,\quad\, j =1,2, \cdots, N_v+1
\end{align*}
where $N_t \Delta t$ reproduces the final time $t_f$ for computation, $[v_1,v_{N_v+1}]$ recovers the velocity domain $[v_{min},v_{max}]$. For space, we assign grid points $x_i$ in different ways:
\begin{itemize}
	\item for periodic boundary condition:
	\begin{align}\label{grid}
		x_i = x_0 + (i-1)\Delta x, i = 1, . . . , N_x.
	\end{align}
	\item for free-flow boundary condition:
	\begin{align}\label{shifted grid}
		x_i = x_0 + (i-\frac{1}{2})\Delta x, i = 1, . . . , N_x.
	\end{align}
\end{itemize}

For the description of numerical methods, we use the following notation. For $s=1,\cdots, L$, we define discrete quantities as follows:
\begin{itemize}
	\item Numerical solutions for $s$-species gas:
	\[
	g_{1,ij}^{s,n+1}\approx g_1^s(t^{n+1},x_i,v_j),\quad g_{2,ij}^{s,n+1}\approx g_2^s(t^{n+1},x_i,v_j).
	\]

	\item Discrete macroscopic moments for $s$-species gas:
	\begin{align}\label{discrete macro}
	n_{s,i}^{n+1} &:= \sum_{j} g_{1,ij}^{s,n+1} \Delta v,\cr
	n_{s,i}^{n+1} u_{s,i}^{n+1} &:= \sum_{j} v_j g_{1,ij}^{s,n+1}  \Delta v,\cr
	\frac{3n_{s,i}^{n+1} K_B T_{i}^{n+1}}{2m_s} &:= \sum_{j} g_{1,ij}^{s,n+1} \frac{\left(v_j- u_{s,i}^{n+1}\right)^2}{2} \Delta v + \sum_{j} \frac{g_{2,ij}^{s,n+1}}{2} \Delta v.
	\end{align}
	\item Collision frequencies
	\begin{align*}
		\nu_{sk,i}^n \approx \nu_{sk}(t^n,x_i)
	\end{align*}
	\item Auxiliary parameters	
	\begin{align}\label{discrete auxiliary}
	\begin{split}	
	a_{sk,i}^{n} = \frac{\lambda_{sk} n_{k,i}^{n} m_k}{\nu_{sk,i}^{n}(m_s + m_k)},\quad b_{sk,i}^{n} = \frac{2a_{sk,i}^{n+1} m_s}{m_s + m_k},\quad
	\gamma_{sk,i}^{n} = \frac{m_i a_{sk,i}^{n+1}}{3}\left( \frac{2m_k}{m_s + m_k} -a_{sk,i}^{n+1}\right),\\
	u_{sk,i}^{n}= (1-a_{sk,i}^{n})u_{s,i}^{n} + a_{sk,i}^{n}u_{k,i}^{n},\qquad\qquad\qquad\qquad\qquad\qquad\\
	T_{sk,i}^{n}= (1-b_{sk,i}^{n})T_{s,i}^{n} + b_{sk,i}^{n}T_{k,i}^{n} + \frac{\gamma_{sk,i}^{n}}{K_B}\left(u_{s,i}^{n}-u_{k,i}^{n}\right)^2.\qquad\qquad\qquad\quad
	\end{split}
	\end{align}
\end{itemize}

\subsection{Derivation of a first order SL scheme.}
Semi-Lagrangian schemes for \eqref{chu bgk bbgsp} can be derived by treating the equation in the Lagrangian formation:
\begin{align}
	\frac{dg_p^s}{dt}&= \frac{1}{\varepsilon}  \nu_{ss} \left(n_{s}M_{ss,p}-g_{p}^{s}\right) + \frac{1}{\kappa} \sum_{k\neq s}^{L} \nu_{sk} \left(n_{s}M_{sk,p}-g_{p}^{s}\right),\quad g_p^s(0,x,v)=g_{p,0}^{s}(\tilde{x},v)\\
	\frac{dx}{dt}&=v,\quad x(0)=\tilde{x}\label{chu bgk bbgsp lag}
\end{align}
where $g_{p,0}^{s}(x,v)$ is the initial condition and $M_{sk,p}$ are given in \eqref{chu bgk bbgsp}. Here the analytic solution to \eqref{chu bgk bbgsp lag} is a characteristic curve passing through the position $x$ with velocity $v$ at time $t$. Tracing back along characteristic curve, we reach the location $\tilde{x}$ and we call this {\it characteristic foot}.

A first order implicit scheme can be obtained by applying implicit Euler method to \eqref{chu bgk bbgsp lag}: 
\begin{align}\label{first order g1}
	g_{1,ij}^{s,n+1} = \tilde{g}_{1,ij}^{s,n} + \frac{\Delta t}{\varepsilon}\nu_{ss,i}^{n+1} \left(n_{s,i}^{n+1}M_{ss,1,ij}^{n+1}-g_{1,ij}^{s,n+1}\right) +   \frac{\Delta t}{\kappa}\sum_{k\neq s}^{L} \nu_{sk,i}^{n+1} \left(n_{s,i}^{n+1}M_{sk,1,ij}^{n+1}-g_{1,ij}^{s,n+1}\right),\\\label{first order g2}
	g_{2,ij}^{s,n+1} = \tilde{g}_{2,ij}^{s,n} + \frac{\Delta t }{\varepsilon}\nu_{ss,i}^{n+1} \left(n_{s,i}^{n+1}M_{ss,2,ij}^{n+1}-g_{2,ij}^{s,n+1}\right) +  \frac{\Delta t}{\kappa} \sum_{k\neq s}^{L} \nu_{sk,i}^{n+1} \left(n_{s,i}^{n+1}M_{sk,2,ij}^{n+1}-g_{2,ij}^{s,n+1}\right).
\end{align}
From the first equation, one can obtain
\begin{align}\label{n approx}
	n_{s,i}^{n+1} =\sum_j \tilde{g}_{1,ij}^{s,n+1} \Delta v =: \tilde{n}_{s,i}^{n+1} ,
\end{align}
Using this, we get
\begin{align}\label{nu approx}
	\nu_{sk,i}^{n+1}=\lambda_{sk} n_{k,i}^{n+1}=\lambda_{sk} \tilde{n}_{k,i}^{n}=\tilde{\nu}_{sk,i}^{n},
\end{align}
and
\begin{align}\label{abgamma approx}
	\begin{split}	
		&a_{sk,i}^{n+1} = \frac{\lambda_{sk} n_{k,i}^{n+1} m_k}{\nu_{sk,i}^{n+1}(m_s + m_k)} = \frac{\lambda_{sk} \tilde{n}_{k,i}^{n} m_k}{\tilde{\nu}_{sk,i}^{n}(m_s + m_k)}=:\tilde{a}_{sk,i}^{n},\cr 
		&b_{sk,i}^{n+1} = \frac{2a_{sk,i}^{n+1} m_s}{m_s + m_k} = \frac{2\tilde{a}_{sk,i}^{n+1} m_s}{m_s + m_k}=:\tilde{b}_{sk,i}^{n},\cr
		&\gamma_{sk,i}^{n+1} = \frac{m_s a_{sk,i}^{n+1}}{3}\left( \frac{2m_k}{m_s + m_k} -a_{sk,i}^{n+1}\right) = \frac{m_s \tilde{a}_{sk,i}^{n}}{3}\left( \frac{2m_k}{m_s + m_k} -\tilde{a}_{sk,i}^{n}\right) =: \tilde{\gamma}_{sk,i}^{n}
	\end{split}
\end{align}
To compute $u_{sk,ij}^{n+1}$, we consider
\begin{align*}
n_{s,i}^{n+1}u_{s,i}^{n+1}
&=\sum_j v_jg_{1,ij}^{s,n+1} \Delta v\cr
&= \sum_j v_j \left(\tilde{g}_{1,ij}^{s,n} + \frac{\Delta t}{\varepsilon}\nu_{ss,i}^{n+1}  \left(n_{s,i}^{n+1}M_{ss,1,ij}^{n+1}-g_{1,ij}^{n+1}\right) + \frac{\Delta t }{\kappa} \sum_{k\neq s}^{L} \nu_{sk,i}^{n+1} \left(n_{s,i}^{n+1}M_{sk,1,ij}^{n+1}-g_{1,ij}^{n+1}\right)\right) \Delta v,
\end{align*}
which can be rewritten as
\begin{align*}
n_{s,i}^{n+1}u_{s,i}^{n+1} &= \tilde{n}_{s,i}^{n}\tilde{u}_{s,i}^{n} + \frac{\Delta t}{\kappa} \sum_{k\neq s}^{L} \nu_{sk,i}^{n+1} \left(n_{s,i}^{n+1}u_{sk,i}^{n+1}-n_{s,i}^{n+1}u_{s,i}^{n+1}\right).
\end{align*}
This, together with \eqref{n approx} and \eqref{nu approx}, gives
\begin{align*}
u_{s,i}^{n+1} &= \tilde{u}_{s,i}^{n} + \frac{\Delta t}{\kappa}\sum_{k\neq s}^{L} \tilde{\nu}_{sk,i}^{n} \left(u_{sk,i}^{n+1}-u_{s,i}^{n+1}\right).
\end{align*}
Now, we recall \eqref{discrete auxiliary} to get 
\begin{align}\label{usk}
	u_{sk,i}^{n+1}= (1-a_{sk,i}^{n+1})u_{s,i}^{n+1} + a_{sk,i}^{n+1}u_{k,i}^{n+1},
\end{align}
and use this to obtain
\begin{align*}
u_{s,i}^{n+1} &= \tilde{u}_{s,i}^{n} + \frac{\Delta t}{\kappa} \sum_{k\neq s}^{L}\tilde{\nu}_{sk,i}^{n}  \left( (1-a_{sk,i}^{n+1})u_{s,i}^{n+1} + a_{sk,i}^{n+1}u_{k,i}^{n+1}  -u_{s,i}^{n+1}\right)\cr
&= \tilde{u}_{s,i}^{n} - \frac{\Delta t}{\kappa} \sum_{k=1}^{L} \tilde{\nu}_{sk,i}^{n} a_{sk,i}^{n+1}\left( u_{s,i}^{n+1} - u_{k,i}^{n+1}\right).
\end{align*}
This combined with the relation \eqref{abgamma approx} gives
\begin{align}\label{U compute}
	u_{s,i}^{n+1} + \frac{\Delta t}{\kappa} \sum_{k\neq s}^{L}  \tilde{\nu}_{sk,i}^{n} \tilde{a}_{sk,i}^{n}\left( u_{s,i}^{n+1} - u_{k,i}^{n+1}\right)&= \tilde{u}_{s,i}^{n}.
\end{align}
which can be rewritten as 
\begin{align}\label{Au}
	\underline{A}_i^{n} \underline{u}_i^{n+1} = \underline{\tilde{u}}_i^n,
\end{align}
where the $(s,k)$ entry of the $s\times s$ matrix $\underline{A}_i^{n}$ is given by
\begin{align*}
	\underline{A}_{i,sk}^n=\begin{cases}
		1+ \frac{\Delta t}{\kappa}\sum_{k\neq s}^{L}  \tilde{\nu}_{sk,i}^{n} \tilde{a}_{sk,i}^{n} , \quad \ \text{if} \quad s=k\\
		-\frac{\Delta t}{\kappa} \tilde{\nu}_{sk,i}^{n} \tilde{a}_{sk,i}^{n}, \qquad\qquad\quad \text{if} \quad s\neq k
	\end{cases}
\end{align*}
with
\begin{align*}
	&\underline{u}_i^{n+1}:=(u_{1,i}^{n+1},\,u_{2,i}^{n+1}, \dots,\, u_{L,i}^{n+1})^\top,\quad \underline{\tilde{u}}_i^{n}:=(\tilde{u}_{1,i}^{n},\,\tilde{u}_{2,i}^{n}, \dots,\, \tilde{u}_{L,i}^{n})^\top.
\end{align*}
We note that the matrix $\underline{A}_i^{n}$ is always invertible because we have 
\begin{align*}
	\left|\underline{A}_{i,ss}^n\right| > \sum_{k\neq s}^L \left|\underline{A}_{i,sk}^n\right|,\quad \text{for every $s=1\cdots,L$}.
\end{align*}
Solving \eqref{Au}, we obtain $u_{s,i}^{n+1}$.

To compute $T_{s,i}^{n+1}$ explicitly,
we take a summation \eqref{first order g1} w.r.t. $\left|v_j-u_{s,i}^{n+1}\right|^2/2$ and \eqref{first order g2} w.r.t. $1/2$, respectively. Then, we sum these to obtain
\begin{align}
\frac{3n_{s,i}^{n+1}K_BT_{s,i}^{n+1}}{2m_s} &- \frac{3\tilde{n}_{s,i}^{n}K_B\tilde{T}_{s,i}^{n}}{2m_s} - \frac{\tilde{n}_{s,i}^{n}}{2} \left(u_{s,i}^{n+1} - \tilde{u}_{s,i}^{n}\right)^2 \cr
&= \frac{\Delta t}{\kappa}\sum_{k\neq s}^{L} \nu_{sk,i}^{n+1} \left(\frac{3n_{s,i}^{n+1}K_BT_{sk,i}^{n+1}}{2m_s} + \frac{n_{s,i}^{n+1}}{2} \left(u_{sk,i}^{n+1} - u_{s,i}^{n+1}\right)^2 -\frac{3n_{s,i}^{n+1}K_BT_{s,i}^{n+1}}{2m_s}\right).
\end{align}
Using \eqref{n approx} and \eqref{nu approx}, we get
\begin{align*}
	T_{s,i}^{n+1}&= \tilde{T}_{s,i}^{n} + \frac{m_s}{3K_B} \left(u_{s,i}^{n+1} - \tilde{u}_{s,i}^{n}\right)^2 +  \frac{\Delta t}{\kappa}\sum_{k\neq s}^{L} \tilde{\nu}_{sk,i}^{n} \left(T_{sk,i}^{n+1} + \frac{m_s}{3K_B} \left(u_{sk,i}^{n+1} - u_{s,i}^{n+1}\right)^2 - T_{s,i}^{n+1}\right).
\end{align*}
Now, we use the relation \eqref{discrete auxiliary} to get
$$ T_{sk,i}^{n+1}= (1-b_{sk,i}^{n+1})T_{s,i}^{n+1} + b_{sk,i}^{n+1}T_{k,i}^{n+1} + \frac{\gamma_{sk,i}^{n+1}}{K_B}\left(u_{s,i}^{n+1}-u_{k,i}^{n+1}\right)^2.$$
This, combined with \eqref{usk}, gives
\begin{align*}
T_{s,i}^{n+1}
&= \tilde{T}_{s,i}^{n} + \frac{m_s}{3K_B} \left(u_{s,i}^{n+1} - \tilde{u}_{s,i}^{n}\right)^2 \cr
&+ \frac{\Delta t}{\kappa} \sum_{k\neq s}^{L} \tilde{\nu}_{sk,i}^{n}  \left(b_{sk,i}^{n+1} \left(-T_{s,i}^{n+1} + T_{k,i}^{n+1}\right)  + \left( \frac{\gamma_{sk,i}^{n+1}}{K_B} + \frac{m_s}{3K_B} (a_{sk,i}^{n+1})^2\right)\left(u_{s,i}^{n+1} - u_{k,i}^{n+1}\right)^2\right).
\end{align*}
Using \eqref{abgamma approx}, we can also derive
\begin{align}\label{T compute}
&T_{s,i}^{n+1} + \frac{\Delta t}{\kappa} \sum_{k\neq s}^{L} \tilde{\nu}_{sk,i}^{n} \tilde{b}_{sk,i}^{n} \left(T_{s,i}^{n+1} - T_{k,i}^{n+1}\right)\cr
&= \tilde{T}_{s,i}^{n} + \frac{m_s}{3K_B} \left(u_{s,i}^{n+1} - \tilde{u}_{s,i}^{n}\right)^2 + \frac{\Delta t}{\kappa} \sum_{k\neq s}^{L} \tilde{\nu}_{sk,i}^{n} \left(\frac{\tilde{\gamma}_{sk,i}^{n}}{K_B} + \frac{m_s}{3K_B} (\tilde{a}_{sk,i}^{n})^2\right)\left(u_{s,i}^{n+1} - u_{k,i}^{n+1}\right)^2
\end{align}
This system also reduces to a solvable linear system:  $$\underline{B}_i^{n} \underline{T}_i^{n+1} = \underline{\tilde{\xi}}_i^n,$$
where 
\begin{align*}
	B_{sk}=\begin{cases}
		1+ \frac{\Delta t }{\kappa} \sum_{k\neq s}^{L} \tilde{\nu}_{sk,i}^{n} \tilde{b}_{sk,i}^{n} , \quad \text{if} \quad s=k\\
		-\frac{\Delta t }{\kappa}\tilde{\nu}_{sk,i}^{n}  \tilde{b}_{sk,i}^{n} , \quad \text{if} \quad s\neq k
	\end{cases}
\end{align*}
\begin{align*}
	&\underline{T}_i^{n+1}:=(T_{1,i}^{n+1},\,T_{2,i}^{n+1}, \dots,\, T_{L,i}^{n+1})^\top,\cr
	&\underline{\tilde{\xi}}_i^{n}:=(\tilde{\xi}_{1,i}^{n},\,\tilde{\xi}_{2,i}^{n}, \dots,\, \tilde{\xi}_{L,i}^{n})^\top,\cr
	&\tilde{\xi}_{s,i}^{n}:= \tilde{T}_{s,i}^{n} + \frac{m_s}{3K_B} \left(u_{s,i}^{n+1} - \tilde{u}_{s,i}^{n}\right)^2 + \frac{\Delta t}{\kappa} \sum_{k\neq s}^{L} \tilde{\nu}_{sk,i}^{n} \left(\frac{\tilde{\gamma}_{sk,i}^{n}}{K_B} + \frac{m_s}{3K_B} (\tilde{a}_{sk,i}^{n})^2\right)\left(u_{s,i}^{n+1} - u_{k,i}^{n+1}\right)^2.
\end{align*}

Since \eqref{U compute} and \eqref{T compute} allows one to compute $u_{s,i}^{n+1}$ and $T_{s,i}^{n+1}$ explicitly, it also gives
\[
M_{sk,1,ij}^{n+1} = M(\text{v};u_{sk,i}^{n+1},\frac{K_BT_{sk,i}^{n+1}}{m_s}),\quad M_{sk,2,ij}^{n+1} =2\frac{K_BT_{sk,i}^{n+1}}{m_s} M_{sk,1,ij}^{n+1}
\]
Then, our numerical solution reads
\begin{align}\label{first order scheme}
	g_{p,ij}^{s,n+1} &= \frac{\tilde{g}_{p,ij}^{s,n} + \frac{\Delta t}{\varepsilon}\nu_{ss,i}^{n+1} n_{s,i}^{n+1}M_{ss,p,ij}^{n+1} +   \frac{\Delta t}{\kappa}\sum_{k\neq s}^{L} \nu_{sk,i}^{n+1} n_{s,i}^{n+1}M_{sk,p,ij}^{n+1}}{1+\frac{\Delta t}{\varepsilon}\nu_{ss,i}^{n+1} +   \frac{\Delta t}{\kappa}\sum_{k\neq s}^{L} \nu_{sk,i}^{n+1}},
\end{align}
where $p=1,2$.

\subsubsection*{Algorithm of first order SL scheme}
\begin{enumerate}
	\item Reconstruct $\tilde{g}_{p,ij}^{s,n}$ at $x_i-v_j\Delta t$ from $\{g_{p,ij}^{s,n}\}$..
	
	\item Compute $\tilde{n}_i^{n}$, $\tilde{u}_i^{n}$, $\tilde{T}_{s,i}^{n}$, $\tilde{\nu}_{sk,i}^{n}$, $\tilde{a}_i^{n}$, $\tilde{b}_i^{n}$, $\tilde{\gamma}_i^{n}$ from $\{\tilde{g}_{p,ij}^{s,n}\}$, 
	
	\item Compute $n_i^{n+1}$, $u_i^{n+1}$, $T_{s,i}^{n+1}$, $\nu_{sk,i}^{n+1}$, $a_i^{n+1}$, $b_i^{n+1}$, $\gamma_i^{n+1}$ with $\tilde{n}_i^{n}$, $\tilde{u}_i^{n}$, $\tilde{T}_{s,i}^{n}$, $\tilde{\nu}_{sk,i}^{n}$, $\tilde{a}_i^{n}$, $\tilde{b}_i^{n}$, $\tilde{\gamma}_i^{n}$ 
	\item Construct $M_{sk,p,ij}^{n+1}$ 
	\item Update solutions with \eqref{first order scheme}.
\end{enumerate}

\begin{remark}
	When we reconstruct solutions on the characteristic feet, the linear interpolation will be sufficient to achieve first order accuracy.
\end{remark}


\section{High order semi-Lagrangian methods}\label{sec high order SL}

\subsection{Time discretization}
In this section, we present high order SL schemes for \eqref{chu bgk bbgsp} by applying L-stable diagonally implicit Runge-Kutta methods (DIRK) or backward difference formula methods (BDF) with high order spatial reconstructions. High order methods can be derived as in the first order scheme, so we only present its algorithm.

\subsubsection{DIRK methods}
DIRK methods can be represented by Butcher's table \cite{HW}:
\begin{align*}
	\begin{array}{c|c}
		c\> &\> A \\
		\hline\\[-3mm]
		&\> b^{\top} 
	\end{array}
\end{align*}
where $A$ is a $s \times s$ lower triangle matrix whose diagonal entries are non-zero constant, $c= (c_1,...,c_s)^{\top}$ and
$b=(b_1,...,b_s)^{\top}$ are coefficients vectors with $c_i = \sum_{j = 1}^s a_{ij}$, for $i = 1,…, s$. In this paper, we use for numerical simulations the following two DIRK methods adopted in \cite{GRS}:
\begin{align}\label{Butcher}
	\text{DIRK2}=\begin{array}{c|c c}
		\alpha \>& \alpha & \>\>0 \\
		1 \>& 1-\alpha & \>\>\alpha \\
		\hline
		\>& 1-\alpha & \>\>\alpha
	\end{array}, \quad 
	\text{DIRK3}=\begin{array}{c|c c c}
		\gamma & \gamma & \>\>0 & \>\>0\\[1mm]
		\frac{1+\gamma}{2} & \frac{1-\gamma}{2} & \>\>\gamma & \>\>0\\[1mm]
		1 & 1-\delta-\gamma & \>\>\delta & \>\>\gamma\\
		\hline
		& 1-\delta-\gamma & \>\>\delta & \>\>\gamma
	\end{array}
\end{align}
where $\alpha= 1-\frac{  \sqrt{2}}{2}$, $\gamma=0.4358665215$ and $\delta =-0.644363171$. For simplicity, we introduce the following notations:
\begin{itemize}
	\item Consider a backward-characteristic which comes from $x_i$ with $v_j$ at time $t^n + c_m \Delta t$. We call this $m$th characteristic. Then, the location of its characteristic at the $\ell$th stage is given by
	:
	\[
	x_{i,j}^{m, \ell}:= x_{i}-(c_m-c_\ell) v_{j}\Delta t
	\]
	Note that $x_{i,j}^{m, m}=x_i$.
	\item  For $p=1,2$, the $\ell$th stage value of the distribution function along the $m$th backward-characteristic which comes from $x_i$ with $v_j$ at time $t^n + c_k \Delta t$ is given by :
	\begin{align*}
		g_{p,ij}^{s,(m,\ell)} :=   g_p^s(x_{i,j}^{m, \ell}, v_{j}, t^n + c_\ell\Delta t).
	\end{align*}
	\item For $p=1,2$, let us define $K_p^s$ by 
	\begin{align*}
	K_p^s:= \frac{1}{\varepsilon} \nu_{ss} \left(n_{s}M_{ss,p}-g_{p}^{s}\right) + \frac{1}{\kappa} \sum_{k\neq s}^{L} \nu_{sk} \left(n_{s}M_{sk,p}-g_{p}^{s}\right).
	\end{align*}
	
	\item For $p=1,2$, the $\ell$-th stage value of $K_p^s$ along the $m$-th backward-characteristic which comes from $x_i$ with $v_j$ at time $t^n + c_k \Delta t$:	
	\begin{align*}
	K_{p,ij}^{s,(m,\ell)} :=   K_p^s(x_{i,j}^{m, \ell}, v_{j}, t^n + c_\ell\Delta t).
	\end{align*}

	
%
\end{itemize}

Using the same argument of the first order SL scheme, all comoutations for 
$\nu$-stage DIRK based SL methods
can be performed explicitly. For $p=1,2$, let us consider $m$-stage values of $g_p^s$:
\begin{align}\label{high order dirk}
	\begin{split}
		g_{p,ij}^{s,(m,m)} &= \tilde{g}_{p,ij}^{s,(m,0)} + \Delta t\sum_{\ell=1}^{m-1} a_{m\ell}K_{p,ij}^{s,(m,\ell)}\cr
		&+ \frac{\Delta t }{\varepsilon}a_{mm} \nu_{ss,ij}^{(m,m)} \left(n_{s,i}^{(m,m)}M_{ss,p,ij}^{(m,m)}-g_{1,ij}^{s,(m,m)}\right)\cr
		&+ \frac{\Delta t }{\kappa}\sum_{k\neq s}^{L}a_{mm} \nu_{sk,ij}^{(m,m)} \left(n_{s,i}^{(m,m)}M_{sk,p,ij}^{(m,m)}-g_{1,ij}^{s,(m,m)}\right)
	\end{split}
\end{align}



In case of $\nu=1$, the implicit Euler method has $a_{11}=1$, and one can update solution with the scheme \eqref{first order scheme}.

In case of $\nu>1$, the computation of $m$-stage value $g_{p,ij}^{s,(m,m)}$, can be understood as repetition of the first order method. However, one should be careful in that the term, $\tilde{g}_{p,ij}^{s,(m,0)} + \Delta t\sum_{\ell=1}^{m-1} a_{m\ell}K_{p,ij}^{s,(m,\ell)}$, now plays the of $\tilde{g}_{p,ij}^{s,n}$ and the coefficient of $\Delta t$ is multiplied by $a_{mm}$. We also note that reconstructions of $g$ and $K$ on characteristic feet will be necessary.
To sum up, DIRK based methods can be updated as follows:

\subsubsection*{Algorithm of $\nu$-stage DIRK based SL scheme}
\begin{enumerate}
	\item Reconstruct $\tilde{g}_{p,ij}^{s,(m,0)}$ at $x_i-v_jc_m\Delta t$ from $\{g_{p,ij}^{s,n}\}$.
	\item Compute $\tilde{g}_{p,ij}^{s,m,*}:= \tilde{g}_{p,ij}^{s,(m,0)} + \Delta t\sum_{\ell=1}^{m-1} a_{m\ell}K_{p,ij}^{s,(m,\ell)}$.
	
	\item Compute $\tilde{n}_i^{m,*}$, $\tilde{u}_i^{m,*}$, $\tilde{T}_{s,i}^{m,*}$, $\tilde{\nu}_{sk,i}^{m,*}$, $\tilde{a}_i^{m,*}$, $\tilde{b}_i^{m,*}$, $\tilde{\gamma}_i^{m,*}$ from $\{\tilde{g}_{p,ij}^{s,m,*}\}$
	
	\item Compute $n_i^{(m,m)}$, $u_i^{(m,m)}$, $T_{s,i}^{(m,m)}$, $\nu_{sk,i}^{(m,m)}$, $a_i^{(m,m)}$, $b_i^{(m,m)}$, $\gamma_i^{(m,m)}$  with $\tilde{n}_i^{m,*}$, $\tilde{u}_i^{m,*}$, $\tilde{T}_{s,i}^{m,*}$, $\tilde{\nu}_{sk,i}^{m,*}$, $\tilde{a}_i^{m,*}$, $\tilde{b}_i^{m,*}$, $\tilde{\gamma}_i^{n}$ 
	
	\item Construct $M_{sk,p,ij}^{(m,m)}$
	\item Update solutions with \eqref{high order dirk}.
\end{enumerate}

\begin{remark}$~$
	\begin{enumerate}
		\item For high order spatial reconstruction, we adopt a conservative reconstruction presented in \cite{BCRY2,BCRY3}.
		\item We note that the relation \eqref{high order dirk} is linear because we explicitly precompute $M_{sk,p,ij}^{(m,m)}$ before computing $g_{p,ij}^{s,(m,m)}$. 
		\item We compute numerical solutions by the last stage value  $g_{p,ij}^{s,n+1}=g_{p,ij}^{s,(\nu,\nu)}$ becasue we consider stiffly accurate DIRK schemes.
	\end{enumerate}
\end{remark}
%
%


\subsubsection{BDF methods}
An efficient class of high order SL  schemes can be designed with linear multi-step methods such as Backward difference formula (BDF) in \cite{HW}. 
The BDF method of order $s$ is represented by
\begin{align}\label{BDFscheme}
	\begin{array}{l}
		\displaystyle	BDF: y^{n+1} = \sum_{k=1}^{s} \alpha_k y^{n+1-k} +  \beta_s \Delta t g(y^{n+1},t_{n+1}).
	\end{array}
\end{align} 
where $\alpha_k$ and $\beta_s$ are constants depending on $s$. In particular, we are to use second and third order methods:
\begin{align}\label{BDFschemes}
	\begin{array}{l}
		\displaystyle	BDF2: \alpha_1=\frac{4}{3},~ \alpha_2=-\frac{1}{3},~ \beta_2=\frac{2}{3}
		,\\[3mm]
		\displaystyle	BDF3:
		\alpha_1=\frac{18}{11},~ \alpha_2=-\frac{9}{11},~ \alpha_3=\frac{2}{11},~ \beta_3=\frac{6}{11}
	\end{array}
\end{align} 

For BDF based methods, we use the following notation:
\begin{itemize}
	\item In case of BDF methods, we consider a backward-characteristic which comes from $x_i$ with $v_j$ at time $t^{n+1}$. Along this characteristic, its characteristic foot at time $t^{n+1-k}$ will be located at
	$x_{i,j}^k:=x_i-k v_j \Delta{t}$ for $k=1,\cdots,s$.
	\item  For $p=1,2$, the value of distribution functions located on $x_{i,j}^k$ at time $t^{n+1-k}$ along the backward-characteristic corresponding to $x_i$ with $v_j$ will be denoted by
	\begin{align*}
		g_{p,ij}^{s,n+1-k} :=   g_p^s(x_{i,j}^k, v_{j}, t^{n+1-k}).
	\end{align*}
	\item For $p=1,2$, the $\ell$-th stage value of $K_p^s$ along the $m$-th backward-characteristic which comes from $x_i$ with $v_j$ at time $t^n + c_k \Delta t$:	
	\begin{align*}
		K_{p,ij}^{s,n+1-k} :=   K_p^s(x_{i,j}^k, v_{j}, t^{n+1-k}).
	\end{align*}	
\end{itemize}



With $s$-order BDF methods, the SL scheme reads
\begin{align}\label{high order bdf}
	\begin{split}
		g_{p,ij}^{s,n+1} &= \sum_{k=1}^{\nu} \alpha_k \tilde{g}_{p,ij}^{s,n+1-k}+\frac{\Delta t}{\varepsilon}  \beta_\nu\nu_{sk,i}^{n+1} \left(n_{s,i}^{n+1}M_{sk,p,ij}^{n+1}-g_{p,ij}^{s,n+1}\right)\cr 
		&\quad + \frac{\Delta t}{\kappa} \sum_{k\neq s}^{L} \beta_\nu \nu_{sk,i}^{n+1} \left(n_{s,i}^{n+1}M_{sk,p,ij}^{n+1}-g_{p,ij}^{s,n+1}\right).
	\end{split}
\end{align}
\subsubsection*{Algorithm of $s$-order BDF based SL scheme}
\begin{enumerate}
	\item Reconstruct $\tilde{g}_{p,ij}^{s,n+1-k}$ at $x_{i,j}^k:=x_i-kv_j\Delta t$ from $\{g_{p,ij}^{s,n+1-k}\}$.
	\item Compute $\tilde{g}_{p,ij}^{s,*}:= \sum_{k=1}^{s} \alpha_k \tilde{g}_{p,ij}^{s,n+1-k}+\frac{\Delta t}{\varepsilon}  \beta_s\nu_{sk,i}^{n+1} \left(n_{s,i}^{n+1}M_{sk,p,ij}^{n+1}-g_{p,ij}^{s,n+1}\right)$.
	
	\item Compute $\tilde{n}_i^{*}$, $\tilde{u}_i^{*}$, $\tilde{T}_{s,i}^{*}$, $\tilde{\nu}_{sk,i}^{*}$, $\tilde{a}_i^{*}$, $\tilde{b}_i^{*}$, $\tilde{\gamma}_i^{*}$ from $\{\tilde{g}_{p,ij}^{s,*}\}.$
	
	\item Compute $n_i^{n+1}$, $u_i^{n+1}$, $T_{s,i}^{n+1}$, $\nu_{sk,i}^{n+1}$, $a_i^{n+1}$, $b_i^{n+1}$, $\gamma_i^{n+1}$ with $\tilde{n}_i^{*}$, $\tilde{u}_i^{*}$, $\tilde{T}_{s,i}^{*}$, $\tilde{\nu}_{sk,i}^{*}$, $\tilde{a}_i^{*}$, $\tilde{b}_i^{*}$, $\tilde{\gamma}_i^{n}$ 
	
	\item Construct $M_{sk,p,ij}^{n+1}$ 
	\item Update solutions with \eqref{high order bdf}.
\end{enumerate}

\begin{figure}[htbp]
	\centering
	\begin{subfigure}[b]{0.43\linewidth}
		\includegraphics[width=1\linewidth]{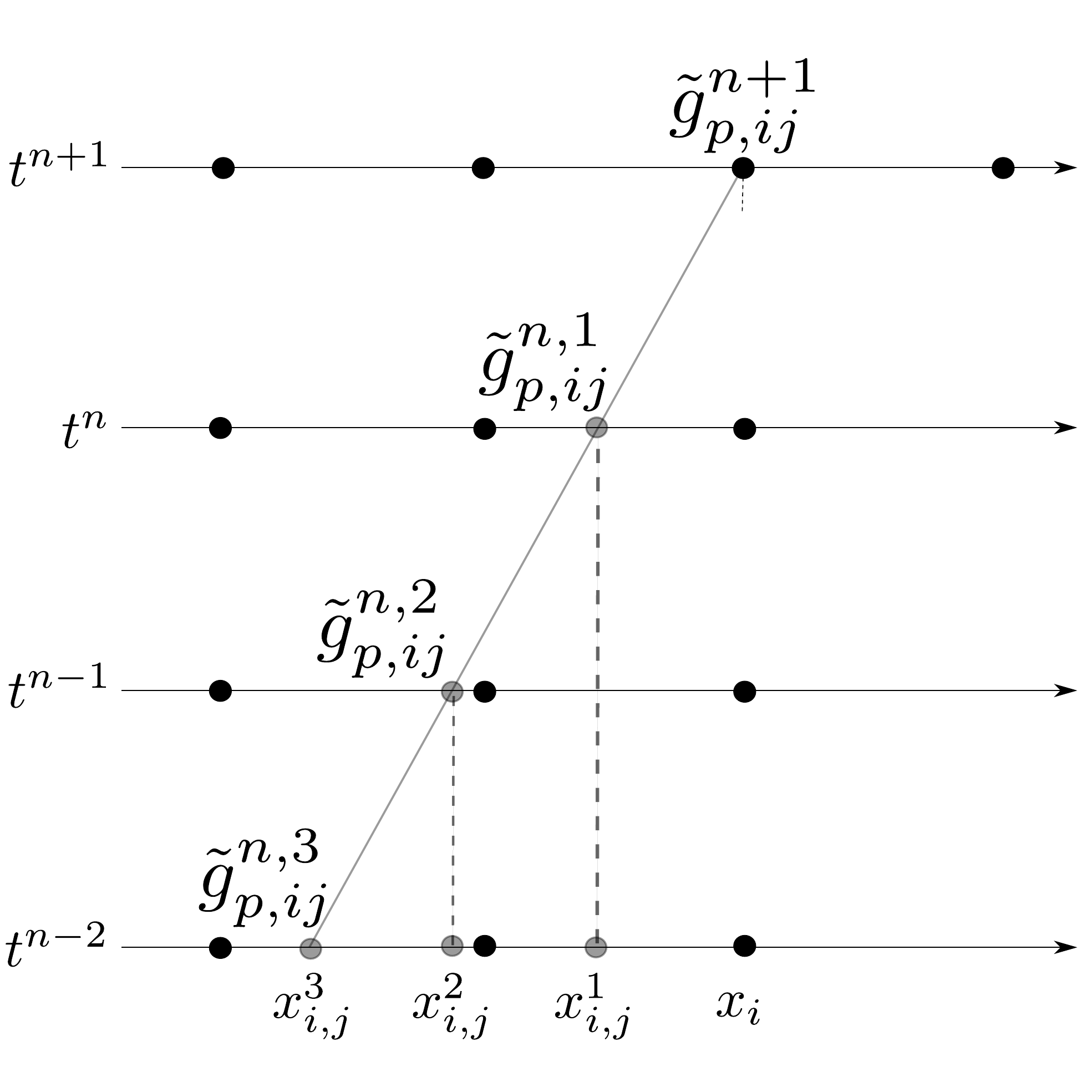}
	\end{subfigure}
	\caption{BDF3 based SL method. On grey points, $\tilde{g}_{p,ij}^{n,k}$ will be reconstructed for each $k=1,2,3$ and $p=1,2$ from $\{g_{p,ij}^{n+1-k}\}$. }\label{fig BDF3}
\end{figure}

In Figure \ref{fig BDF3}, we describe the schematic of BDF3 based SL scheme. Note that a characteristic starting from position $x_i$ with velocity $v_j$ at time $t^{n+1}$ passes through three points $x_{i,j}^k$, $k=1,2,3$ where reconstructions are necessary.

\section{Conservative reconstruction}\label{sec conservative recon}
In this section, we briefly present a conservative reconstruction technique introduced in \cite{BCRY2,BCRY3}. The main idea of the approach is to take a sliding average of a basic reconstruction. As in literatures \cite{BCRY2,BCRY3}, we also take CWENO reconstructions as basic reconstructions. In the following description, we consider two cases $k=2,4$, which correspond to CWENO23 \cite{LPR2} and CWENO35 \cite{CPSV}, respectively.
We summarize the procedure of the conservative reconstruction as follows:
\begin{enumerate}
	\item Given point-wise values $\{u_i\}$ on grid points $\{x_i\}$, compute a CWENO reconstruction for each cell $I_i=[x_i-\Delta x/2, x_i+\Delta x/2] $:
	\begin{align*}
		R_i(x)= \sum_{\ell=0}^k \frac{R_i^{(\ell)}}{\ell !}(x-x_i)^{\ell}.
	\end{align*}
	This yields a CWENO reconstruction on the spatial domain: 
	$$R(x)=\sum_i R_i(x) \chi_i(x),$$ 
	where $\chi_i(x)$ is an indicator function defined on cell $I_i$.
	\item For any $\theta \in [0,1)$, approximate $u(x_i + \theta \Delta x)$ by $Q(x_i+\theta \Delta x)$:
	$$ u(x_i + \theta \Delta x) \approx Q(x_i+\theta \Delta x):=\frac{1}{\Delta x} \int_{x_i + \theta \Delta x-\Delta x/2}^{x_i + \theta \Delta x-\Delta x/2} R(y)dy$$
	whose explicit form is given by 
	\begin{align*}
		Q(x_i+\theta \Delta x)= \sum_{\ell=0}^k (\Delta x)^{\ell} \left(\alpha_{\ell}(\theta)R_{i}^{(\ell)} +\beta_\ell(\theta)R_{i+1}^{(\ell)}\right),
	\end{align*}
	where  
	\begin{align*}
		\alpha_\ell(\theta) = \frac{1- (2\theta -1)^{\ell+1}}{2^{\ell+1}(\ell+1)!}, \quad 
		\beta_\ell(\theta) = \frac{(2\theta -1)^{\ell+1} - (-1)^{\ell+1} }  {2^{\ell+1}(\ell+1)!}.
	\end{align*}
\end{enumerate}
In the following section, the resulting conservative reconstruction based on CWENO23 $(k=2)$ will be denoted by Q-CWENO23. Similarly, we call the conservative reconstruction based on CWENO35 $(k=4)$ by Q-CWENO35.
\section{Numerical tests}\label{sec numerical test}
In this section, we verify the performance of high-order conservative semi-Lagrangian methods to \eqref{chu bgk bbgsp}.
For tests, we use high order schemes dealt with in sections \ref{sec high order SL}-\ref{sec conservative recon}. In particular, we consider the following combination of time integration methods and spatial reconstructions:
\begin{enumerate}
	\item RK2-QCW23: DIRK2 with Q-CWENO23.
	\item RK3-QCW35: DIRK3 with Q-CWENO35.
	\item BDF2-QCW23: BDF2 with Q-CWENO23.
	\item BDF3-QCW35: BDF3 with Q-CWENO35.
\end{enumerate}
Here the Q-CWENO23 and Q-CWENO35 reconstructions are high order conservative reconstruction adopted in \cite{BCRY2,BCRY3}. Throughout this section, we determine a time step $\Delta t$ based on the CFL number defined by
\[
\text{CFL}= \max_j|v_j|\frac{\Delta t}{\Delta x}.
\]
Considering the relation in \eqref{positivity condition}, 
throughout numerical experiments we consider collision frequencies of the following form:
$\nu_{sk}= \lambda_{sk} n_{k}$ which will be divided by $\varepsilon$ or $\kappa$.


\subsection{Accuracy test}\label{sec accuracy}
As in \cite{ADGG,GRS2}, we consider gas mixtures of four monoatomic gases whose molecular masses are given by
\[
m_1= 58.5,\quad  m_2 = 18,\quad  m_3 = 40,\quad  m_4 = 36.5.
\]
The values of $\lambda_{sk}$ are:
\begin{align}\label{lambda values}
	\lambda_{11}= 5,\quad \lambda_{12}= 6,\quad \lambda_{13}=2,\quad \lambda_{14}=7& \cr 
	\lambda_{22}= 4,\quad \lambda_{23}=5,\quad \lambda_{24}=8& \cr 
	\lambda_{33}=4,\quad \lambda_{34}=3& \cr 
	\lambda_{44}=6&,		
\end{align}
where $\lambda_{sk}= \lambda_{ks}$ for $s, k = 1,\dots, 4$. Here we set $\kappa=\varepsilon$. We take as initial data Maxwellians whose macroscopic fields are 
\begin{align}\label{initial acc}
\begin{split}
	n_0^s(x)=\frac{1}{m_s},\quad 
T_0^s(x)=\frac{4}{\sum_{s=1}^4 n_0^s},\qquad\qquad\qquad\qquad\cr
u_0^s(x)= \frac{s}{\sigma_s} \left[   \exp\left(-\left(\sigma_s x- 1 + \frac{s}{3} \right)^2\right) + \exp\left(-\left(\sigma_s x+ 3 - \frac{s}{10} \right)^2\right)\right].
\end{split}
\end{align}
for $s = 1, \cdots, 4$, where $\sigma_s = (10, 13, 16, 19)$. The periodic condition is imposed on the space domain $[-1,1]$ and velocity domain is truncated by $[-15,15]$. We compute numerical solutions using a time step determined by CFL$=2$ up to $t_f=0.2$. Here we use $N_v=60$ and $N_x=40,\, 80,\, 160,\, 320$. For each scheme, we expect to observe both the order of time integration and that of spatial reconstruction.

In Table \ref{tab accuracy}, we report the numerical errors and convergence rates for global number density based on the relative $L^1$-norm. Here we confirm that RK2, BDF2 and BDF3 based schemes attain desired accuracies.  However, we observe the order reduction with the RK3 based scheme in the limit  $\varepsilon \rightarrow 0$, which commonly happens in the stiff problem when we adopt RK methods of order greater than two.

\begin{center}
	\begin{table}[h]
		\small
		\centering
		{\begin{tabular}{|cccccccccc|}
				\hline
				\multicolumn{1}{ |c }{}&
				\multicolumn{1}{ |c| }{}& \multicolumn{2}{ c|  }{$\varepsilon=10^{-5}$} & \multicolumn{2}{ c|  }{$\varepsilon=10^{-4}$} & \multicolumn{2}{ c|  }{$\varepsilon=10^{-3}$} & \multicolumn{2}{ c|  }{$\varepsilon=10^{-2}$}  \\ \hline
				\multicolumn{1}{ |c }{}&
				\multicolumn{1}{ |c|  }{$N_x$,$2N_x$} &
				\multicolumn{1}{ c  }{error} &
				\multicolumn{1}{ c|  }{rate} &
				\multicolumn{1}{ c  }{error} &
				\multicolumn{1}{ c|  }{rate} &
				\multicolumn{1}{ c  }{error} &
				\multicolumn{1}{ c|  }{rate}&
				\multicolumn{1}{ c  }{error} &
				\multicolumn{1}{ c|  }{rate}   \\ 
				\hline
				\hline	
				\multicolumn{1}{ |c }{RK2}&
				\multicolumn{1}{ |c|  }{$40,80$}&3.01e-03
				&1.97
				&2.95e-03
				&2.03
				&2.50e-03
				&2.35
				&8.88e-04
				&2.72
				\\
				\multicolumn{1}{ |c }{QW23}&
				\multicolumn{1}{ |c|  }{$80,160$}&7.69e-04
				&2.07
				&7.22e-04
				&2.15
				&4.90e-04
				&2.65
				&1.35e-04
				&2.93
				\\
				\multicolumn{1}{ |c }{}&
				\multicolumn{1}{ |c|  }{$160,320$}&1.84e-04
				&&1.63e-04
				&&7.81e-05
				&&1.78e-05&
				\\
				\hline
				\hline	
				\multicolumn{1}{ |c }{BDF2}&
				\multicolumn{1}{ |c|  }{$40,80$}&3.58e-03
				&1.97
				&3.54e-03
				&2.01
				&3.04e-03
				&2.28
				&1.15e-03
				&2.59
				\\
				\multicolumn{1}{ |c }{QW23}&
				\multicolumn{1}{ |c|  }{$80,160$}&9.14e-04
				&1.83
				&8.76e-04
				&1.88
				&6.27e-04
				&2.48
				&1.91e-04
				&2.76
				\\
				\multicolumn{1}{ |c }{}&
				\multicolumn{1}{ |c|  }{$160,320$}&2.57e-04
				&&2.38e-04
				&&1.13e-04
				&&2.82e-05&
				\\
				\hline
				\hline	
				\multicolumn{1}{ |c }{RK3}&
				\multicolumn{1}{ |c|  }{$40,80$}&2.46e-03
				&1.97
				&2.33e-03
				&2.19
				&1.69e-03
				&2.93
				&1.00e-03
				&4.23
				\\
				\multicolumn{1}{ |c }{QW35}&
				\multicolumn{1}{ |c|  }{$80,160$}&6.30e-04
				&1.88
				&5.09e-04
				&2.89
				&2.22e-04
				&4.39
				&5.35e-05
				&4.73
				\\
				\multicolumn{1}{ |c }{}&
				\multicolumn{1}{ |c|  }{$160,320$}&1.71e-04
				&&6.85e-05
				&&1.06e-05
				&&2.02e-06&
				\\
				\hline
				\hline	
				\multicolumn{1}{ |c }{BDF3}&
				\multicolumn{1}{ |c|  }{$40,80$}&2.74e-03
				&2.01
				&2.65e-03
				&2.06
				&2.16e-03
				&2.81
				&7.86e-04
				&4.49
				\\
				\multicolumn{1}{ |c }{QW35}&
				\multicolumn{1}{ |c|  }{$80,160$}&6.79e-04
				&3.43
				&6.37e-04
				&3.02
				&3.08e-04
				&3.58
				&3.51e-05
				&4.92
				\\
				\multicolumn{1}{ |c }{}&
				\multicolumn{1}{ |c|  }{$160,320$}&6.28e-05
				&&7.86e-05
				&&2.58e-05
				&&1.16e-06&
				\\
				\hline
				\hline
		\end{tabular}}
		\caption{Accuracy test for SL methods associated to initial data \eqref{initial acc}.}\label{tab accuracy}
		
	\end{table}
\end{center}

\subsection{Indifferentiability principle}\label{Indifferentiabiliy}
The aim of this section is to check if the scheme maintains the indifferentiabiliy principle of \eqref{bgk bbgsp}. For this, we show a mixture of four equal gases behaves like a single gas, independently of $\varepsilon$. In this test, we assume $\kappa=\varepsilon$.


\subsubsection{Single gas and a mixture of four gases}\label{sec Single gas and four gases}
We start from the comparison of single gas and a mixture of four gases. Here we consider following numerical setting:
\begin{enumerate}
	\item Single gas ($L=1$): Consider
	a single monoatomic gas of molecular mass $m_1= 58.5$ with collision frequency $\lambda_{11}= 5$.
	We take the Maxwellian as initial data which reproduce
	\begin{align}\label{initial acc single 1}
		\begin{split}
			n_0(x)=\frac{4}{m_1},\quad 
			T_0(x)=\frac{4}{n_0^1},\qquad\qquad\qquad\qquad\cr
			u_0(x)= \frac{1}{10} \left[   \exp\left(-\left(10 x- 1 + \frac{1}{3} \right)^2\right) -2 \exp\left(-\left(10 x+ 3 - \frac{1}{10} \right)^2\right)\right].
		\end{split}
	\end{align}
	\item Four gases ($L=4$): For fair comparison with the single gas case, we consider
	a mixture of binary gases whose molecular masses are given by
	\[
	m_1= m_2 = m_3 = m_4 = 58.5,
	\]
	with $\lambda_{sk}= 5$ for all $1\leq s,k\leq L$. We also set initial data by the Maxwellian whose macroscopic fields are given as follows:
	\begin{align}\label{initial acc single 2}
		\begin{split}
			n_0^s(x)=\frac{1}{m_1},\quad 		T_0^s(x)=\frac{4}{\sum_{s=1}^4 n_0^s},\qquad\qquad\qquad\qquad\qquad\cr
			u_0^s(x)= \frac{1}{10} \left[   \exp\left(-\left(10 x- 1 + \frac{1}{3} \right)^2\right) -2 \exp\left(-\left(10 x+ 3 - \frac{1}{10} \right)^2\right)\right].
		\end{split}
	\end{align}
for $s=1,2,3,4$.
\end{enumerate}

For both cases, We impose the periodic boundary condition on the space domain $[-1,1]$ and truncate velocity domain by $[-15,15]$. We compute numerical solutions with $N_v=60$, $N_x=200$ up to a final time $t_f=0.2$. For a time interval $t\in (0,0.02]$, due to the use of not well-prepared initial data we set CFL$=0.2$. After then we continue the computation with CFL$=2$.

\begin{figure}[htbp]
	\centering
	\begin{subfigure}[b]{0.43\linewidth}
		\includegraphics[width=1\linewidth]{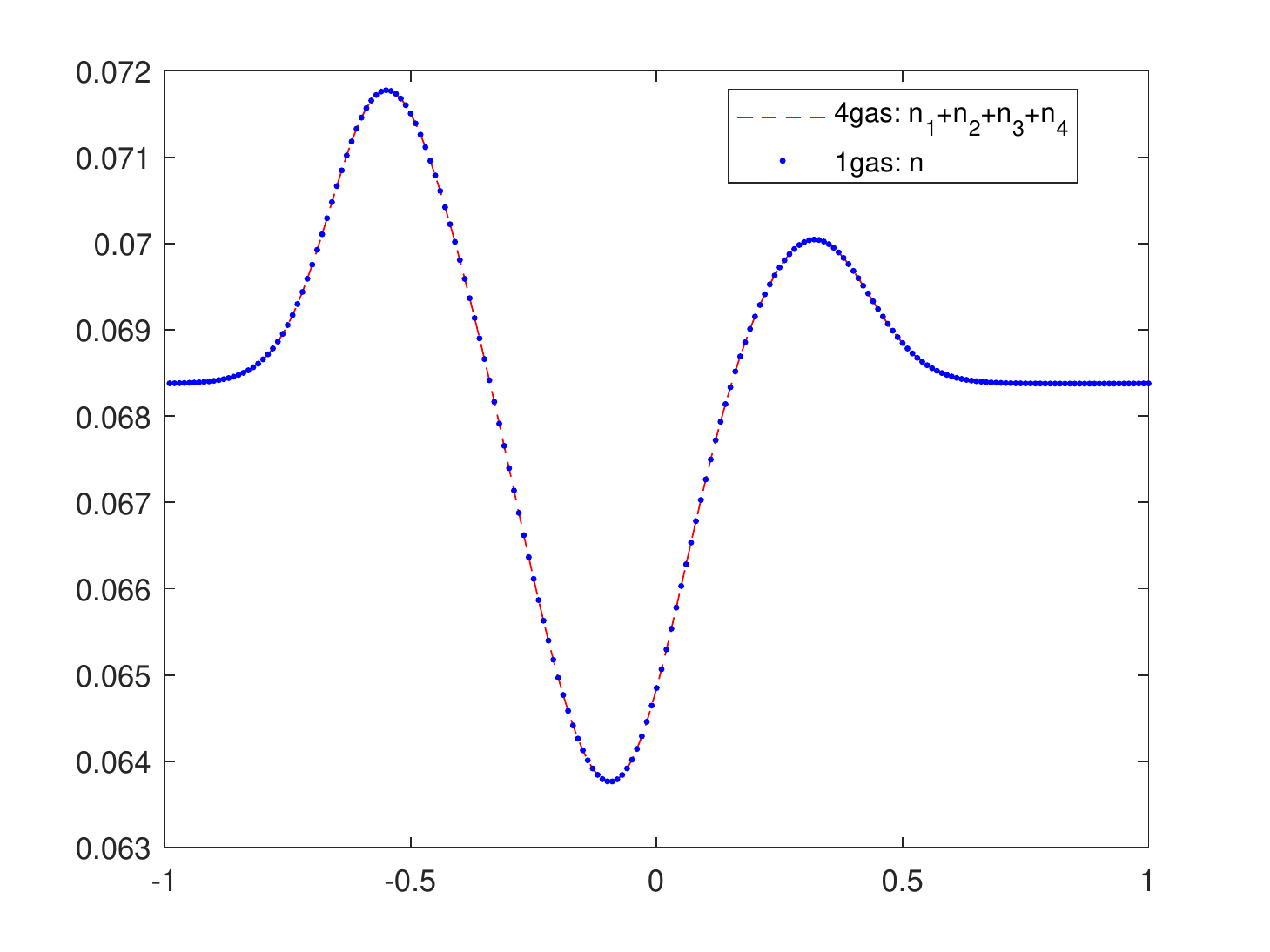}
		\subcaption{Number density}
	\end{subfigure}
	\begin{subfigure}[b]{0.43\linewidth}
		\includegraphics[width=1\linewidth]{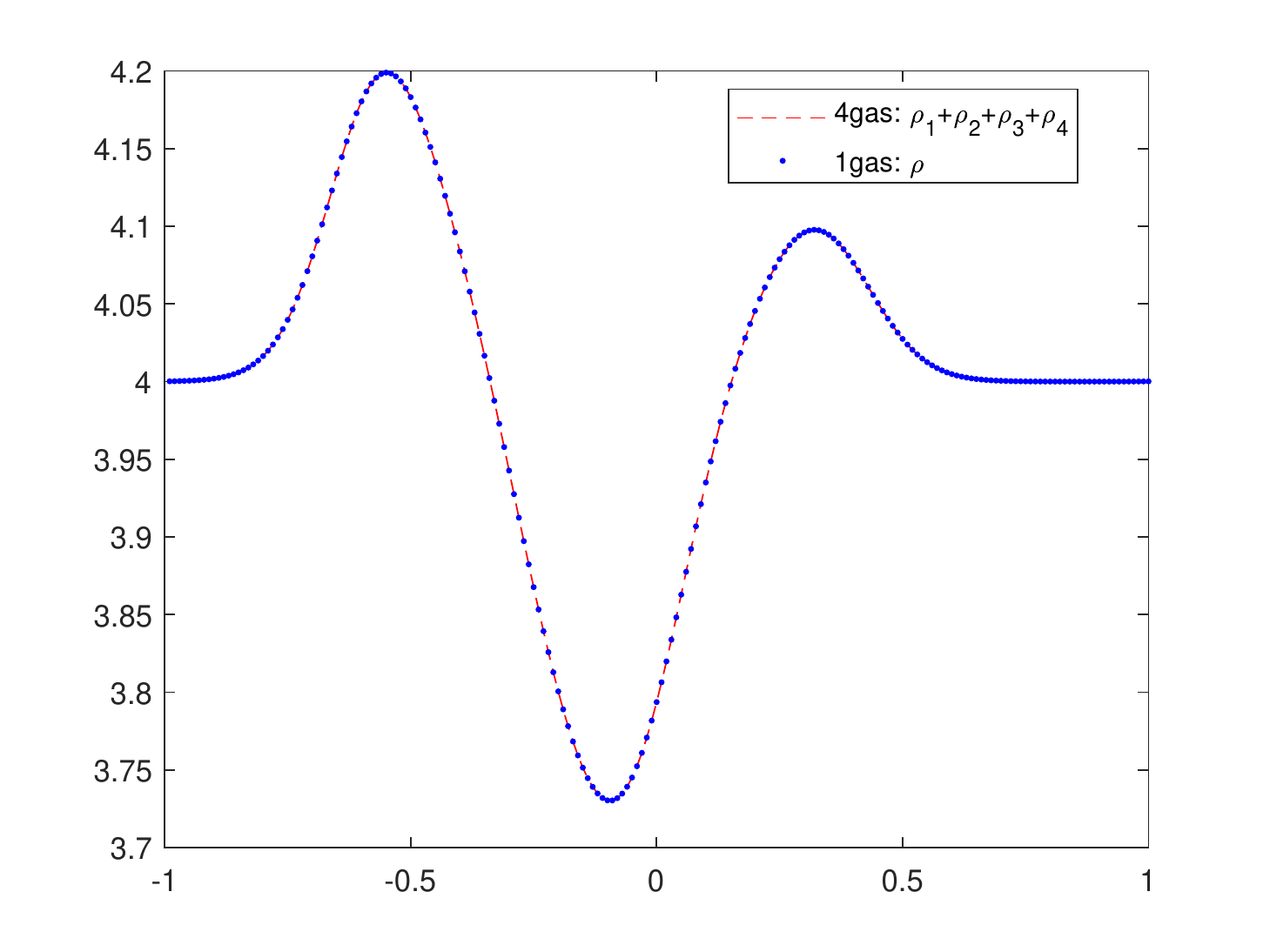}
		\subcaption{Density}
	\end{subfigure}
	\begin{subfigure}[b]{0.43\linewidth}
		\includegraphics[width=1\linewidth]{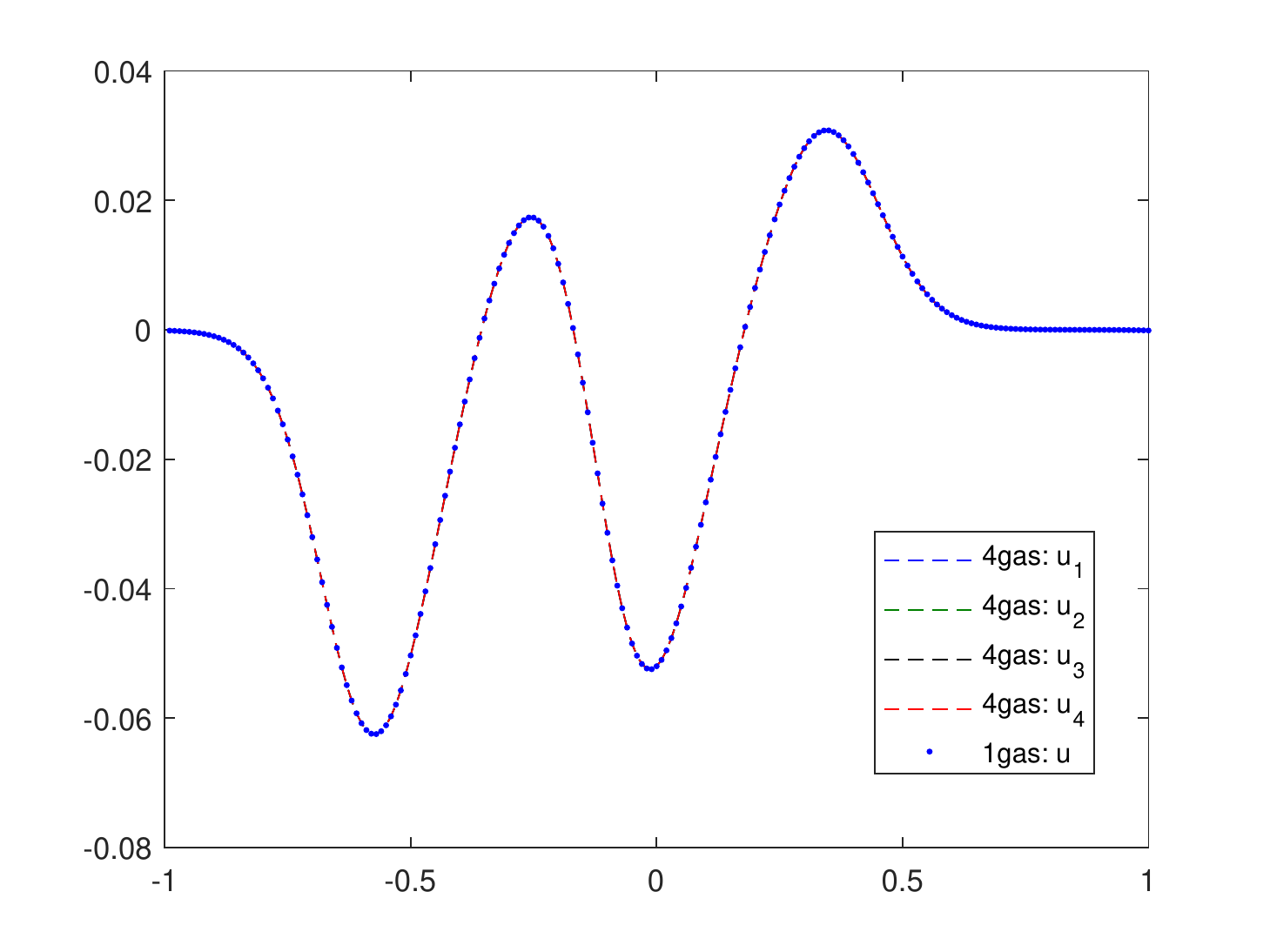}
		\subcaption{Velocity}
	\end{subfigure}			
	\begin{subfigure}[b]{0.43\linewidth}
		\includegraphics[width=1\linewidth]{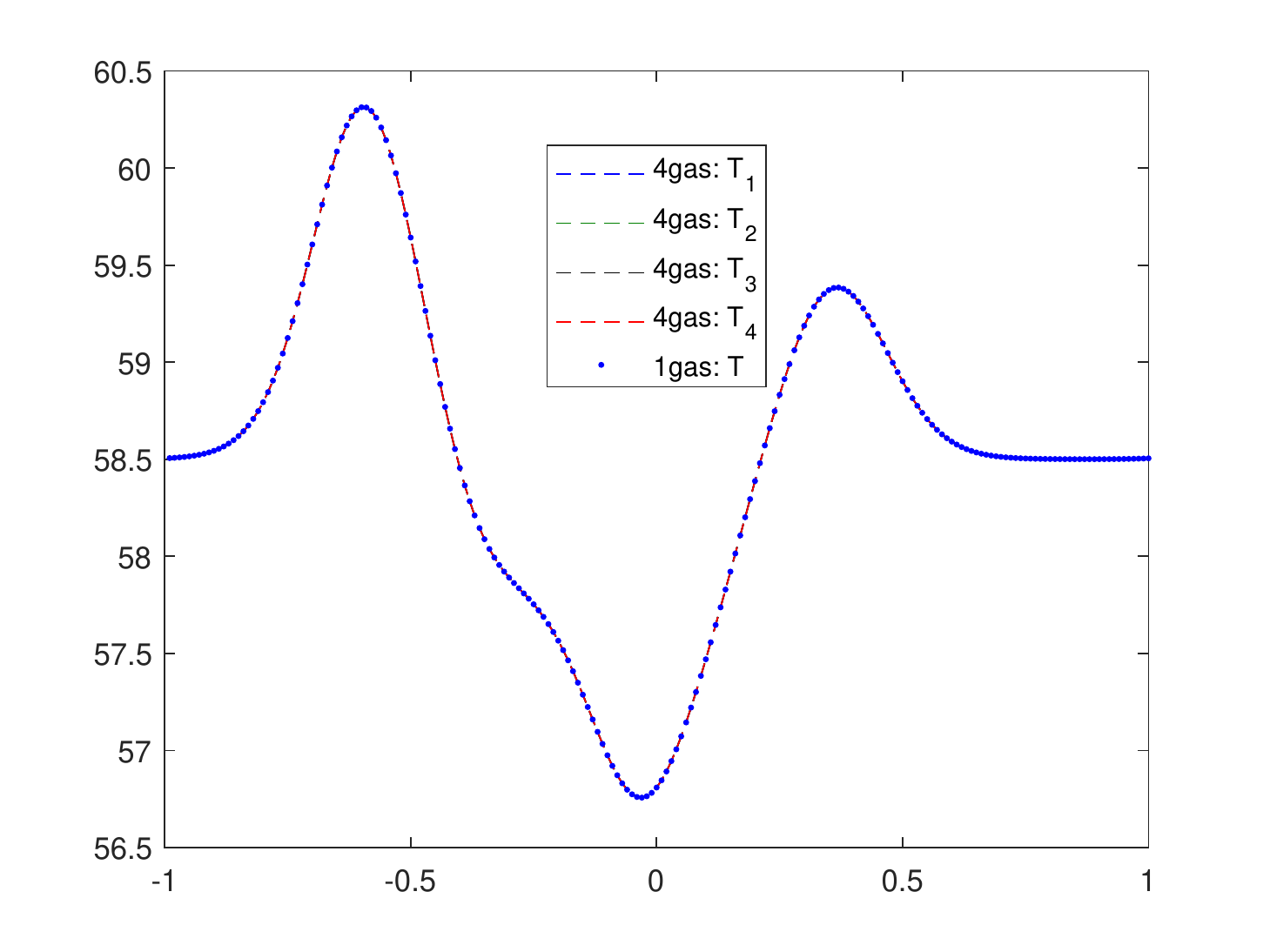}
		\subcaption{Temperature}
	\end{subfigure}			
	\caption{BDF3-QCWENO35 for $\varepsilon=10^{-2}$ 
		associated to initial data \eqref{initial acc single 1} and \eqref{initial acc single 2}.}\label{fig comparison single 01}
\end{figure}

In Tab. \ref{tab id principle error}, we report the 
discrepancy of solutions of single gas and four gases using the relative $L^1$-norm. For various ranges of $\varepsilon=10^{-q}$, $q=2,3,4,5$, we observe that the spatial errors are dominant with respect to time errors and the accuracy of BDF3-QCW35 scheme is confirmed.
In Fig. \ref{fig comparison single 01}, we note that both solutions show good agreements to each other for $q=2$. Similar results are obtained with the other schemes proposed in the paper.

\begin{center}
	\begin{table}[h]
		\centering
		{\begin{tabular}{|cccccccccc|}
				\hline
				\multicolumn{1}{ |c }{}&
				\multicolumn{1}{ |c| }{$N_x$}& \multicolumn{2}{ c|  }{$\varepsilon=10^{-5}$} & \multicolumn{2}{ c|  }{$\varepsilon=10^{-4}$} & \multicolumn{2}{ c|  }{$\varepsilon=10^{-3}$} & \multicolumn{2}{ c|  }{$\varepsilon=10^{-2}$}  \\ 
				\hline
				\multicolumn{1}{ |c }{}&
				\multicolumn{1}{ |c|  }{} &
				\multicolumn{1}{ c|  }{error} &
				\multicolumn{1}{ c|  }{rate} &
				\multicolumn{1}{ c|  }{error} &
				\multicolumn{1}{ c|  }{rate} &
				\multicolumn{1}{ c|  }{error} &
				\multicolumn{1}{ c|  }{rate} &
				\multicolumn{1}{ c|  }{error} &
				\multicolumn{1}{ c|  }{rate}   \\ 		
				\hline
				\hline	
				\multicolumn{1}{ |c }{}&
				\multicolumn{1}{ |c|  }{100}&
				\multicolumn{1}{ c|  }{4.73e-06}&
				\multicolumn{1}{ c|  }{3.26}&
				\multicolumn{1}{ c|  }{4.22e-06}&
				\multicolumn{1}{ c|  }{3.30}&
				\multicolumn{1}{ c|  }{1.93e-06}&
				\multicolumn{1}{ c|  }{3.44}&
				\multicolumn{1}{ c|  }{4.37e-07}&
				\multicolumn{1}{ c|  }{3.80}
				\\
				\multicolumn{1}{ |c }{$n$}&
				\multicolumn{1}{ |c|  }{200}&
				\multicolumn{1}{ c|  }{4.95e-07}&
				\multicolumn{1}{ c|  }{4.17}&
				\multicolumn{1}{ c|  }{4.28e-07}&
				\multicolumn{1}{ c|  }{4.33}&
				\multicolumn{1}{ c|  }{1.78e-07}&
				\multicolumn{1}{ c|  }{4.60}&
				\multicolumn{1}{ c|  }{3.15e-08}&
				\multicolumn{1}{ c|  }{4.86}
				\\
				\multicolumn{1}{ |c }{}&
				\multicolumn{1}{ |c|  }{400}&
				\multicolumn{1}{ c|  }{2.76e-08}&
				\multicolumn{1}{ c|  }{}&
				\multicolumn{1}{ c|  }{2.12e-08}&
				\multicolumn{1}{ c|  }{}&
				\multicolumn{1}{ c|  }{7.36e-09}&
				\multicolumn{1}{ c|  }{}&
				\multicolumn{1}{ c|  }{1.09e-09}&
				\multicolumn{1}{ c|  }{}
				\\
				\hline
				\hline	
				\multicolumn{1}{ |c }{}&
				\multicolumn{1}{ |c|  }{100}&
				\multicolumn{1}{ c|  }{2.68e-04}&
				\multicolumn{1}{ c|  }{3.39}&
				\multicolumn{1}{ c|  }{2.44e-04}&
				\multicolumn{1}{ c|  }{3.47}&
				\multicolumn{1}{ c|  }{1.17e-04}&
				\multicolumn{1}{ c|  }{3.67}&
				\multicolumn{1}{ c|  }{3.23e-05}&
				\multicolumn{1}{ c|  }{4.01}
				\\
				\multicolumn{1}{ |c }{$u$}&
				\multicolumn{1}{ |c|  }{200}&
				\multicolumn{1}{ c|  }{2.55e-05}&
				\multicolumn{1}{ c|  }{4.11}&
				\multicolumn{1}{ c|  }{2.21e-05}&
				\multicolumn{1}{ c|  }{4.28}&
				\multicolumn{1}{ c|  }{9.20e-06}&
				\multicolumn{1}{ c|  }{4.65}&
				\multicolumn{1}{ c|  }{2.01e-06}&
				\multicolumn{1}{ c|  }{4.94}
				\\
				\multicolumn{1}{ |c }{}&
				\multicolumn{1}{ |c|  }{400}&
				\multicolumn{1}{ c|  }{1.48e-06}&
				\multicolumn{1}{ c|  }{}&
				\multicolumn{1}{ c|  }{1.14e-06}&
				\multicolumn{1}{ c|  }{}&
				\multicolumn{1}{ c|  }{3.66e-07}&
				\multicolumn{1}{ c|  }{}&
				\multicolumn{1}{ c|  }{6.55e-08}&
				\multicolumn{1}{ c|  }{}
				\\
				\hline
				\hline	
				\multicolumn{1}{ |c }{}&
				\multicolumn{1}{ |c|  }{100}&
				\multicolumn{1}{ c|  }{3.16e-06}&
				\multicolumn{1}{ c|  }{3.31}&
				\multicolumn{1}{ c|  }{2.83e-06}&
				\multicolumn{1}{ c|  }{3.38}&
				\multicolumn{1}{ c|  }{1.34e-06}&
				\multicolumn{1}{ c|  }{3.55}&
				\multicolumn{1}{ c|  }{2.92e-07}&
				\multicolumn{1}{ c|  }{3.72}
				\\
				\multicolumn{1}{ |c }{$T$}&
				\multicolumn{1}{ |c|  }{200}&
				\multicolumn{1}{ c|  }{3.19e-07}&
				\multicolumn{1}{ c|  }{4.21}&
				\multicolumn{1}{ c|  }{2.71e-07}&
				\multicolumn{1}{ c|  }{4.37}&
				\multicolumn{1}{ c|  }{1.14e-07}&
				\multicolumn{1}{ c|  }{4.61}&
				\multicolumn{1}{ c|  }{2.21e-08}&
				\multicolumn{1}{ c|  }{4.92}
				\\
				\multicolumn{1}{ |c }{}&
				\multicolumn{1}{ |c|  }{400}&
				\multicolumn{1}{ c|  }{1.72e-08}&
				\multicolumn{1}{ c|  }{}&
				\multicolumn{1}{ c|  }{1.31e-08}&
				\multicolumn{1}{ c|  }{}&
				\multicolumn{1}{ c|  }{4.70-09}&
				\multicolumn{1}{ c|  }{}&
				\multicolumn{1}{ c|  }{7.28e-10}&
				\multicolumn{1}{ c|  }{}
				\\
				\hline
				\hline
		\end{tabular}}
		\caption{BDF3-QCW35. The discrepancy of solutions for 4gas and 1gas. Numerical solutions are computed with initial data \eqref{initial acc single 1} and \eqref{initial acc single 2}.}\label{tab id principle error}
	\end{table}
\end{center}

\subsection{Comparison with classical Euler system \eqref{Euler 1d}}\label{sec euler}
In this test, we aim to check if the proposed scheme gives correct solutions in the Euler limit. For this, we consider a test in \cite{GRS2}. The molecular masses and value of $\lambda_{sk}$ follows the choice in \eqref{lambda values}. The Maxwellian is taken as initial data, which reproduces the following macroscopic fields:
\begin{align*}
\left(\rho_0, u_0, p_0\right) =
\begin{cases}
\left(1,\,0,\,5/3\right), \qquad x<0.5\\
\left(1/8,\,0,\,1/6\right), \quad x>0.5,
\end{cases}
\end{align*}
\begin{align*}
\left(\rho_{01}, \rho_{02}, \rho_{03}, \rho_{04}\right) =
\begin{cases}
\left(1/10,\,2/10,\,3/10,\,4/10\right), \quad x<0.5\\
\left(1/80,\,2/80,\,3/80,\,4/80\right), \quad x>0.5.
\end{cases}
\end{align*}
We assume the freeflow boundary condition on the spatial domain $[-1,1]$ and truncate velocity domain by $[-15,15]$. We compute numerical solutions with $N_x=200$ and $N_v=60$ upto $T_f=0.2$. As in the previous test, we set CFL=$0.2$ for $t \in (0,0.2]$ and CFL=$2$ for $t\in (0.2,2]$ for this test. We take $\varepsilon = \kappa = 10^{-6}$

\begin{figure}[htbp]
	\centering
	\begin{subfigure}[b]{0.43\linewidth}
	\includegraphics[width=1\linewidth]{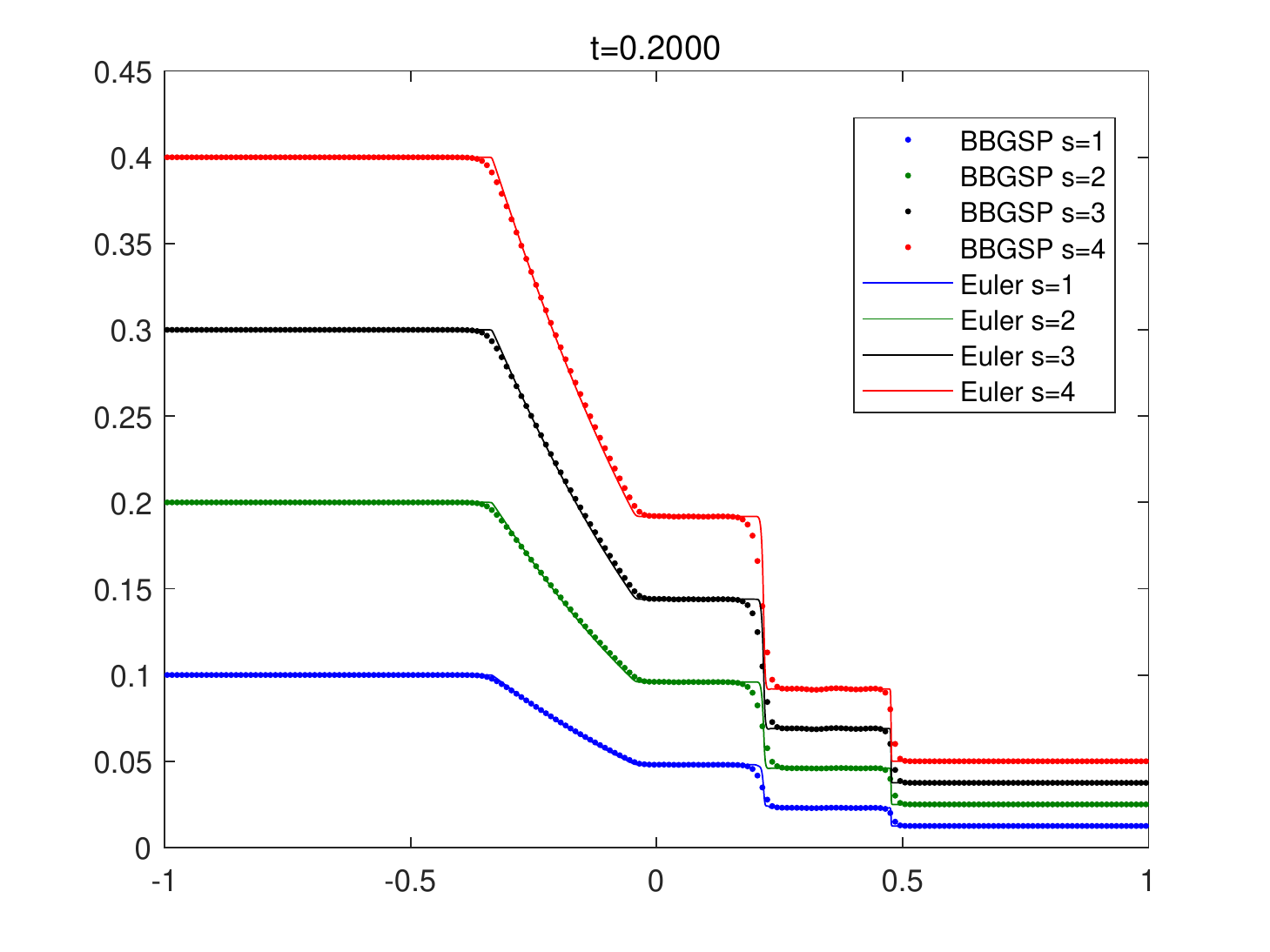}
	\subcaption{Density $\rho_s$, $s=1,2,3,4$}
	\end{subfigure}
	\begin{subfigure}[b]{0.43\linewidth}
		\includegraphics[width=1\linewidth]{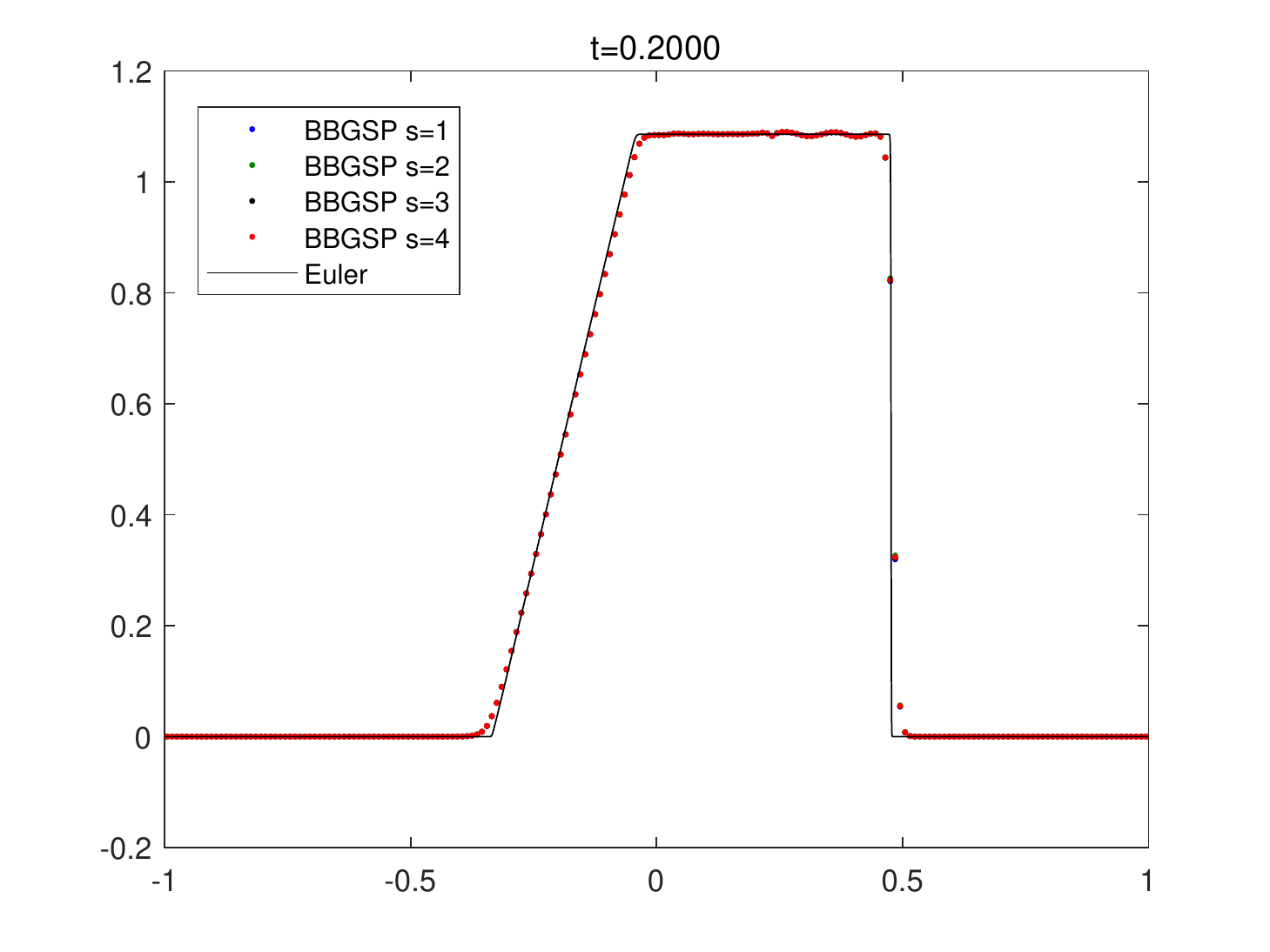}
		\subcaption{Velocity $u_s$, $s=1,2,3,4$ (BGK), $u$ (Euler)}
	\end{subfigure}
	\begin{subfigure}[b]{0.43\linewidth}
		\includegraphics[width=1\linewidth]{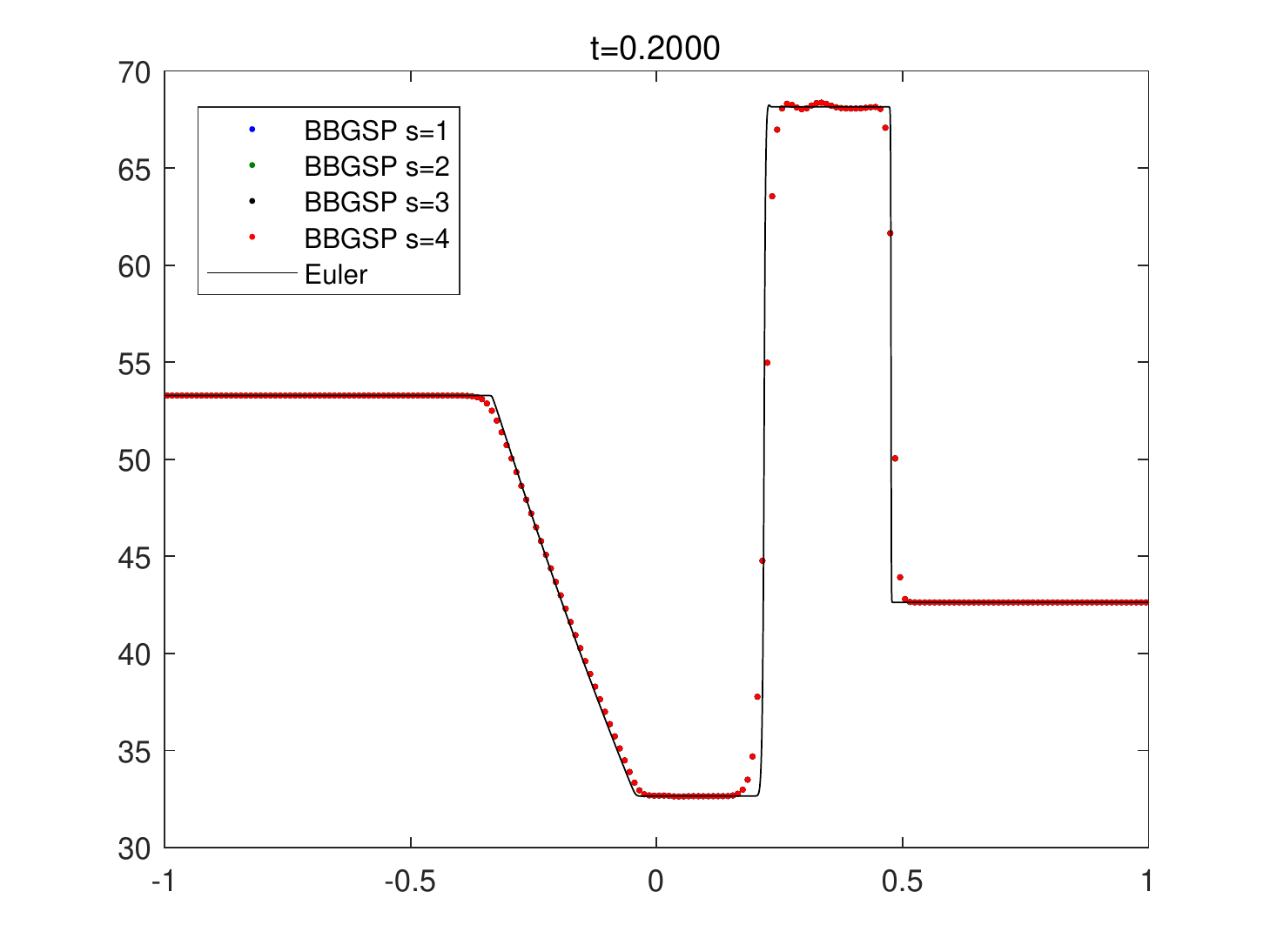}
		\subcaption{Temperature $T_s$, $s=1,2,3,4$ (BGK), $T$ (Euler)}
	\end{subfigure}	
	\begin{subfigure}[b]{0.43\linewidth}
		\includegraphics[width=1\linewidth]{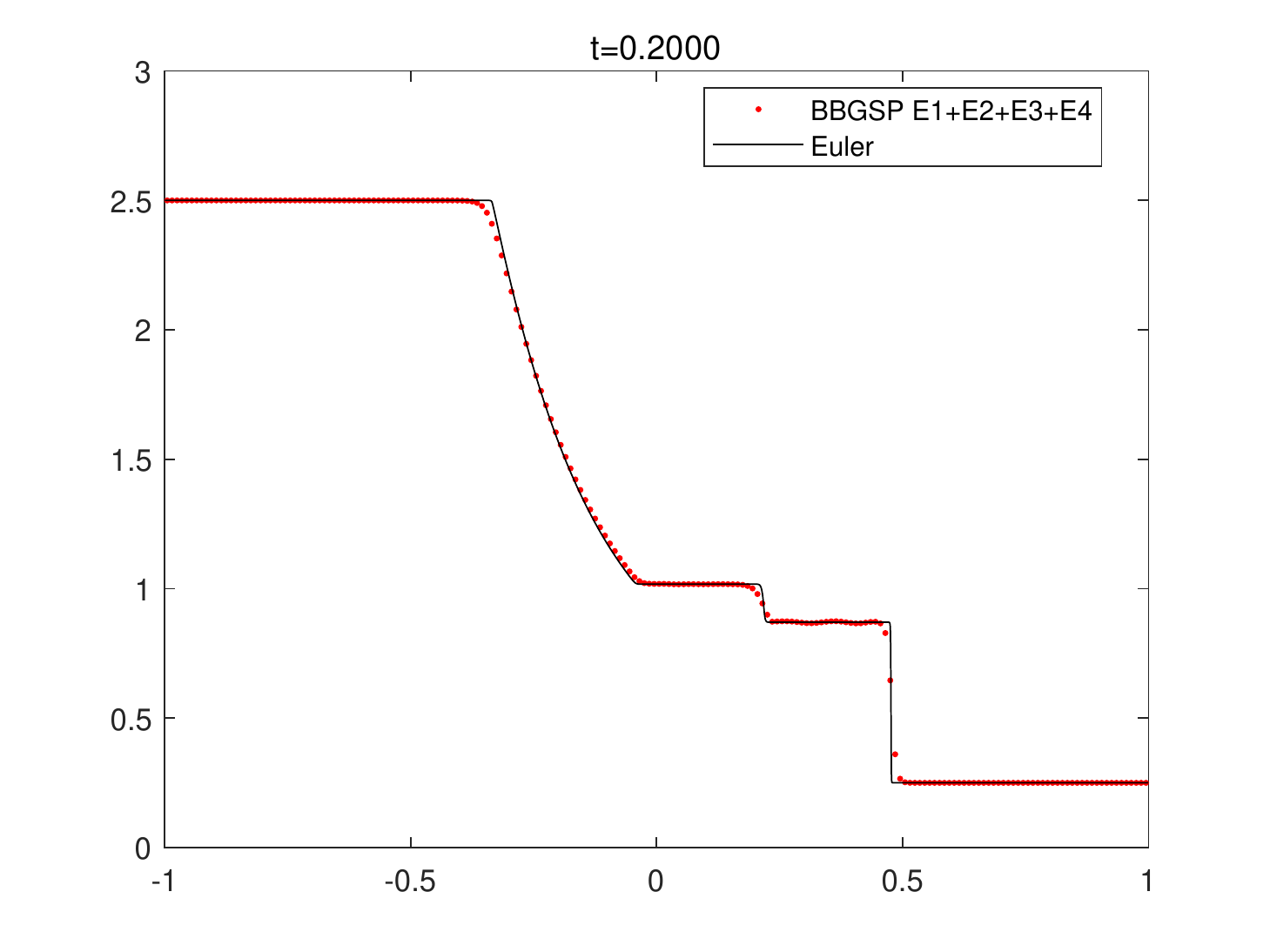}
		\subcaption{Energy $E_1+E_2+E_3+E_4$ (BGK), $E$ (Euler)}
	\end{subfigure}			
	\caption{BDF3-QCW35 for $\varepsilon=10^{-6}$ associated to initial data in section \ref{sec euler}. Reference solutions are plotted by black line $'-'$.}\label{fig euler}
\end{figure}

In Figure \ref{fig euler}, we compare our solution with the reference solution obtained by solving the classical Euler system \eqref{Euler 1d} with explicit RK3 method and WENO23 (with 4000 grid points). The numerical solutions show good agreement with the reference solutions. In particular, we can observe that each species velocity and temperature recovers the same global velocity and temperature of gas mixture. This result also verifies that the proposed scheme has the asymptotic preserving property.


\subsection{Comparison with multi velocities and temperatures Euler system \eqref{Euler multi 1d}}\label{sec euler multi}
In the limit $\kappa \rightarrow 0$, we expect the solution to \eqref{Euler multi 1d} converges to that of \eqref{Euler 1d}. To confirm this asymptotic behavior,
we consider the same numerical setting in section \ref{sec euler}. We set a small fixed value of $\varepsilon=10^{-6}$, and take different values of $\kappa=10^{-q}$, $q=1,2,3,4$. 

\begin{figure}[htbp]
	\centering
	\begin{subfigure}[b]{0.43\linewidth}
		\includegraphics[width=1\linewidth]{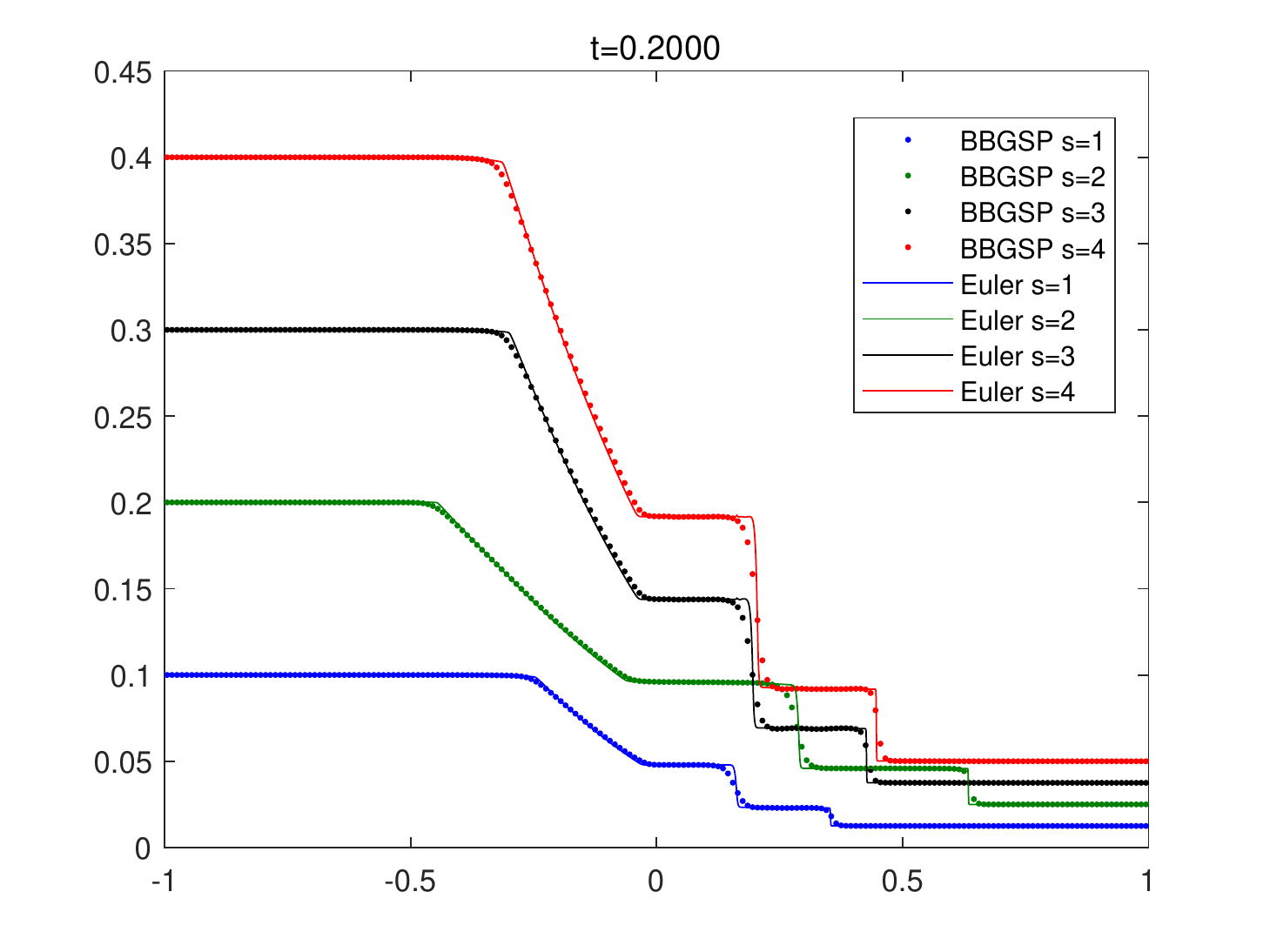}
		\subcaption{Density $\rho_s$, $s=1,2,3,4$}
	\end{subfigure}
	\begin{subfigure}[b]{0.43\linewidth}
		\includegraphics[width=1\linewidth]{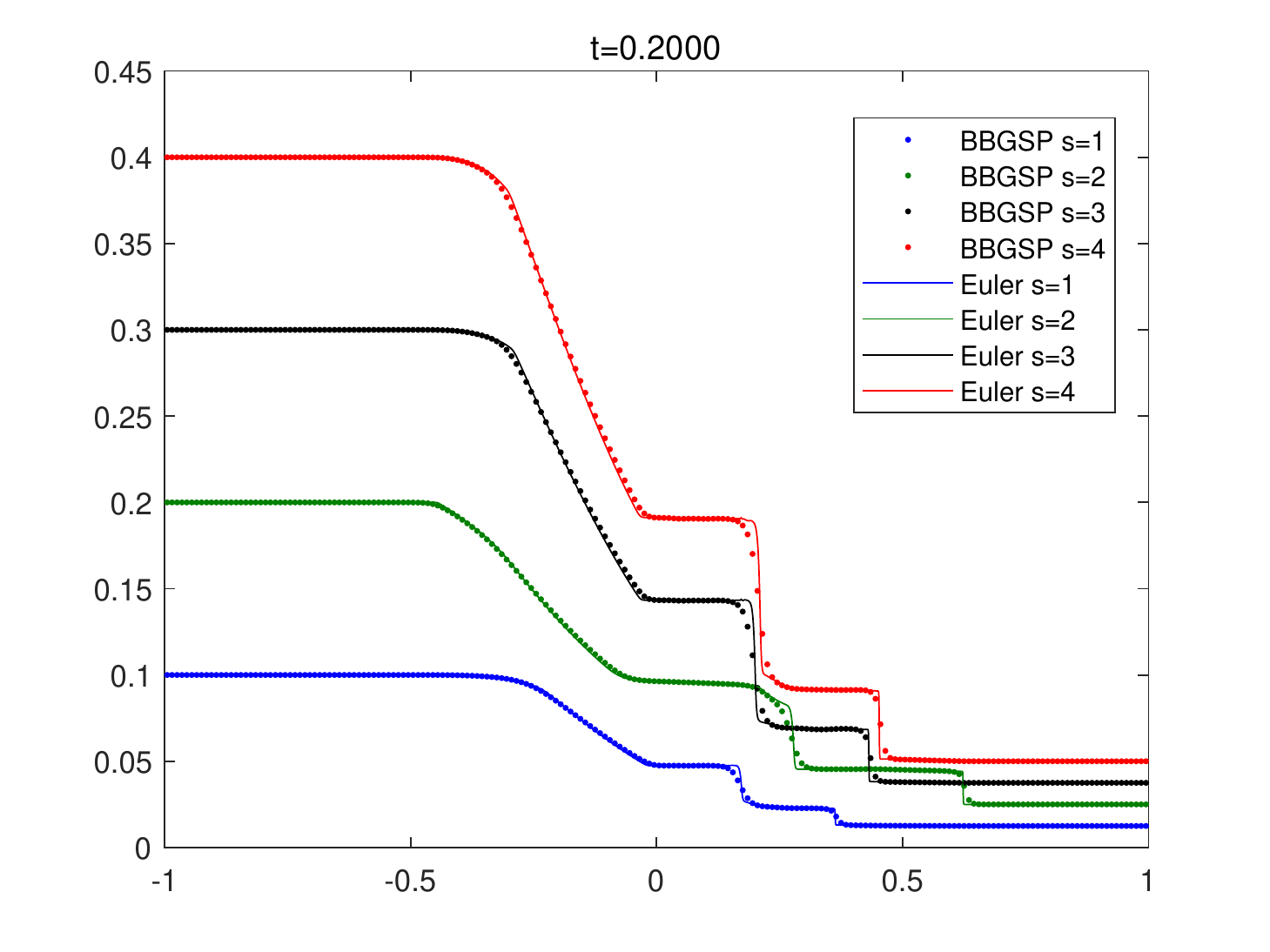}
		\subcaption{Density $\rho_s$, $s=1,2,3,4$}
	\end{subfigure}
\begin{subfigure}[b]{0.43\linewidth}
	\includegraphics[width=1\linewidth]{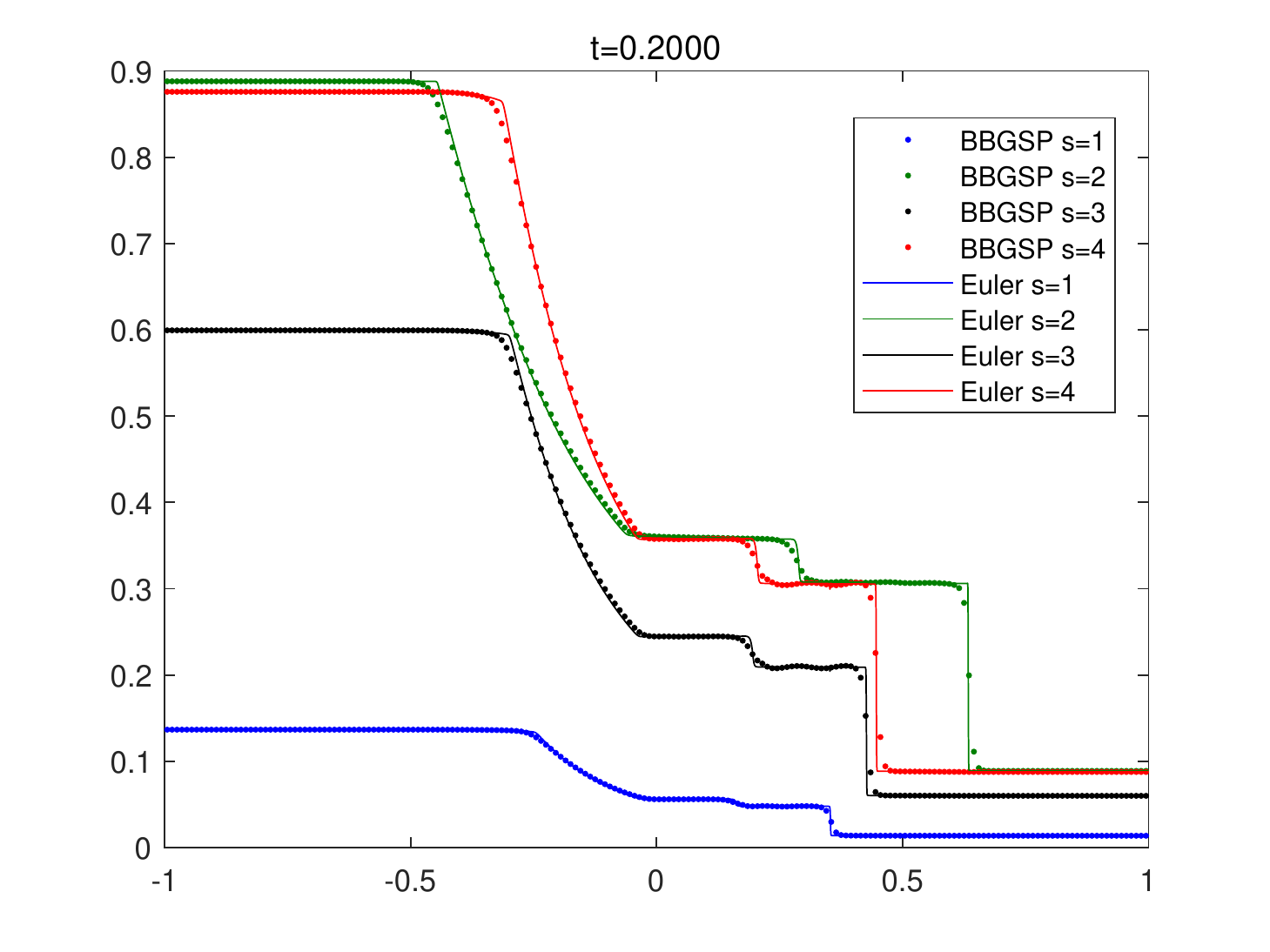}
	\subcaption{Energy $E_s$, $s=1,2,3,4$}
\end{subfigure}
\begin{subfigure}[b]{0.43\linewidth}
	\includegraphics[width=1\linewidth]{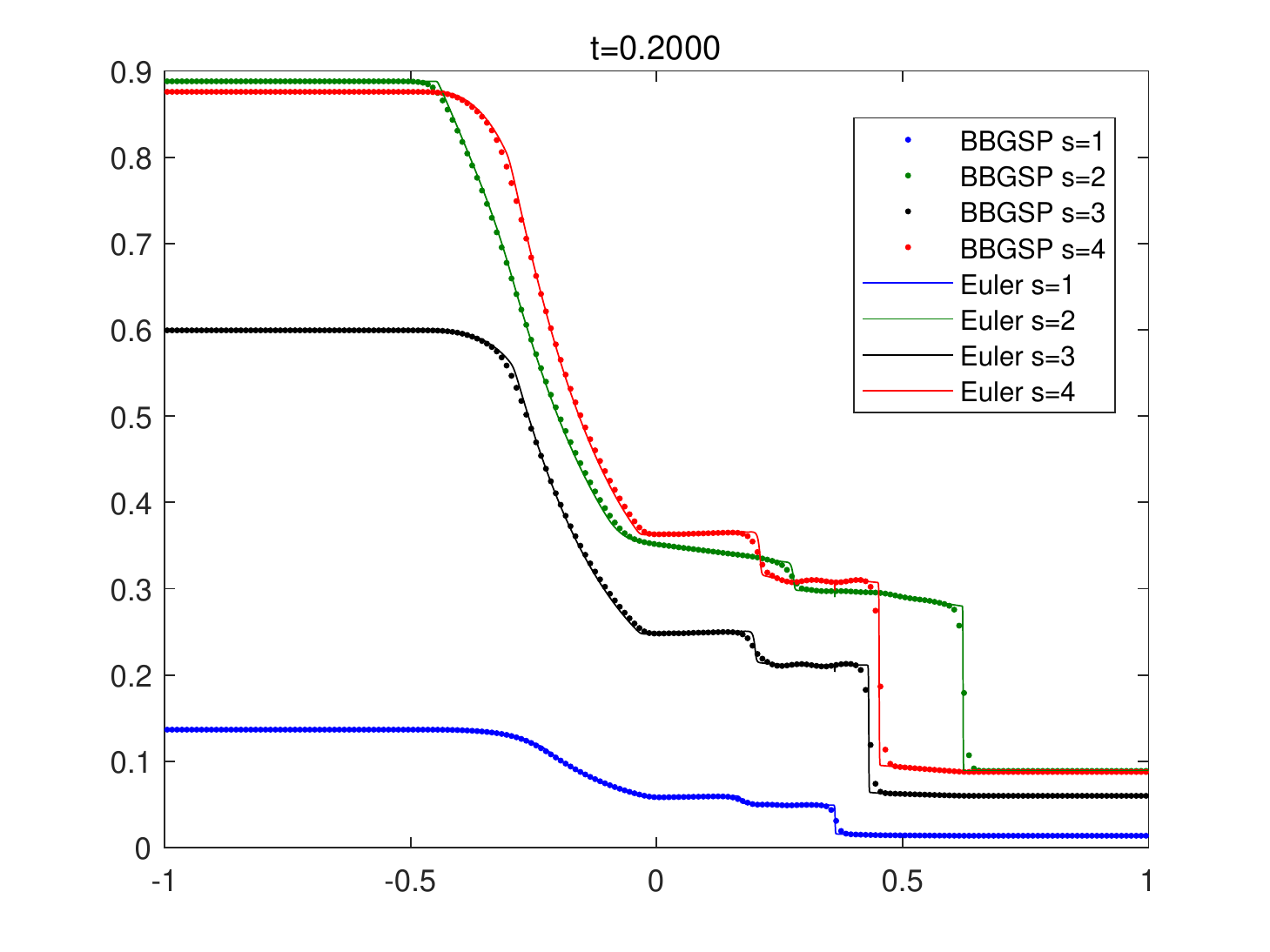}
	\subcaption{Energy $E_s$, $s=1,2,3,4$}
\end{subfigure}
	\begin{subfigure}[b]{0.43\linewidth}
		\includegraphics[width=1\linewidth]{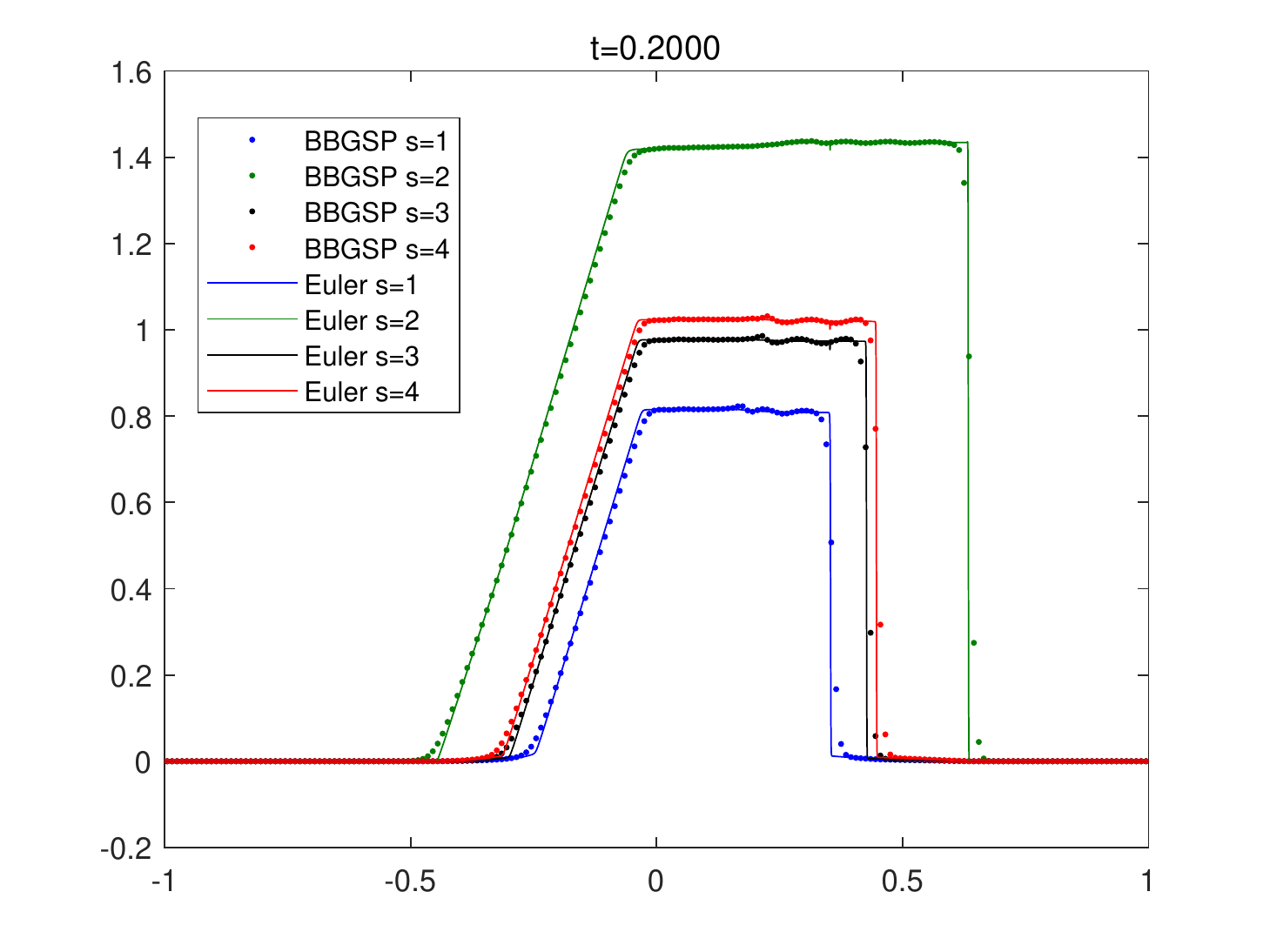}
		\subcaption{Velocity $u_s$, $s=1,2,3,4$}
	\end{subfigure}
	\begin{subfigure}[b]{0.43\linewidth}
		\includegraphics[width=1\linewidth]{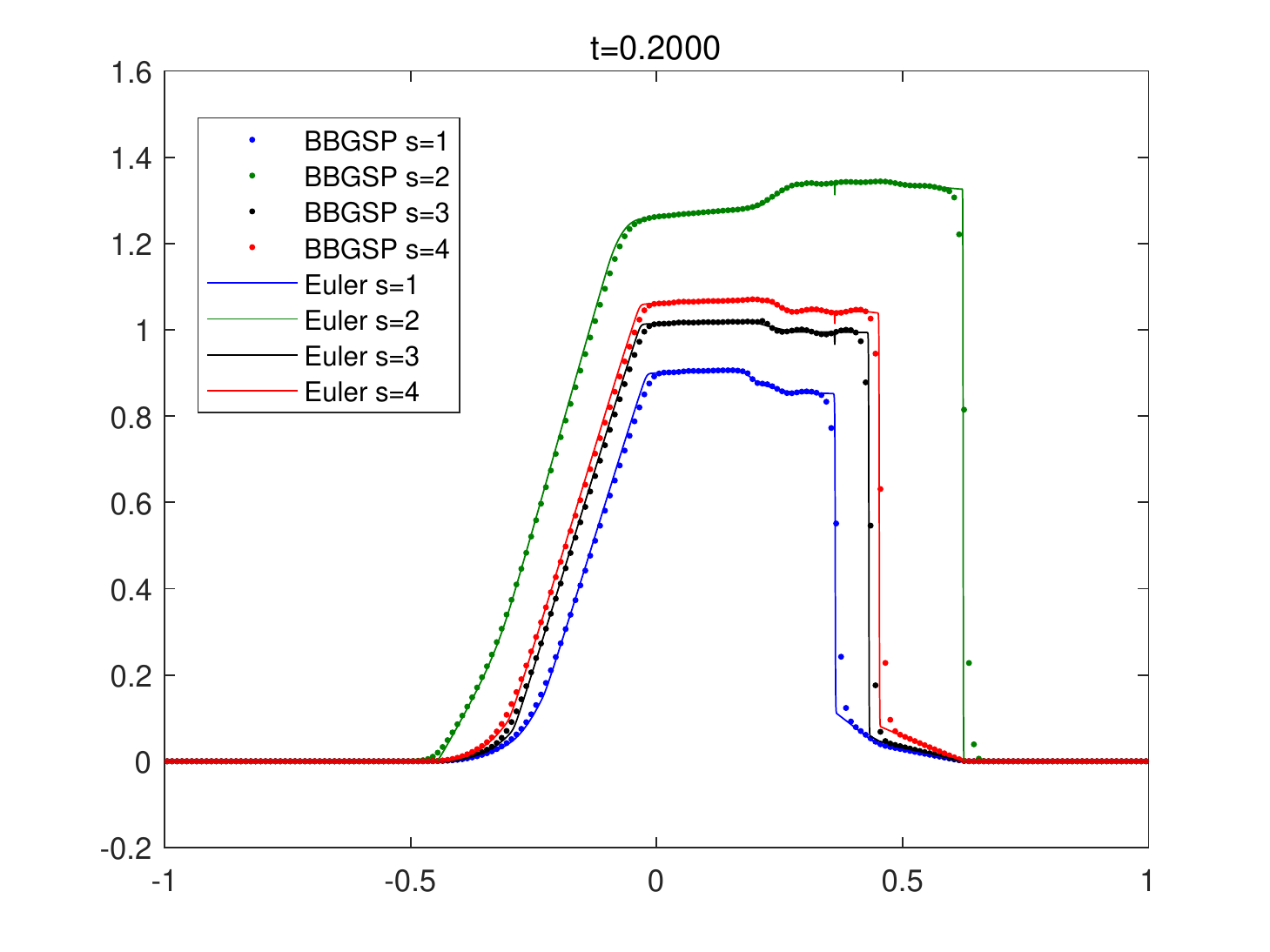}
		\subcaption{Velocity $u_s$, $s=1,2,3,4$}
	\end{subfigure}
	\begin{subfigure}[b]{0.43\linewidth}
		\includegraphics[width=1\linewidth]{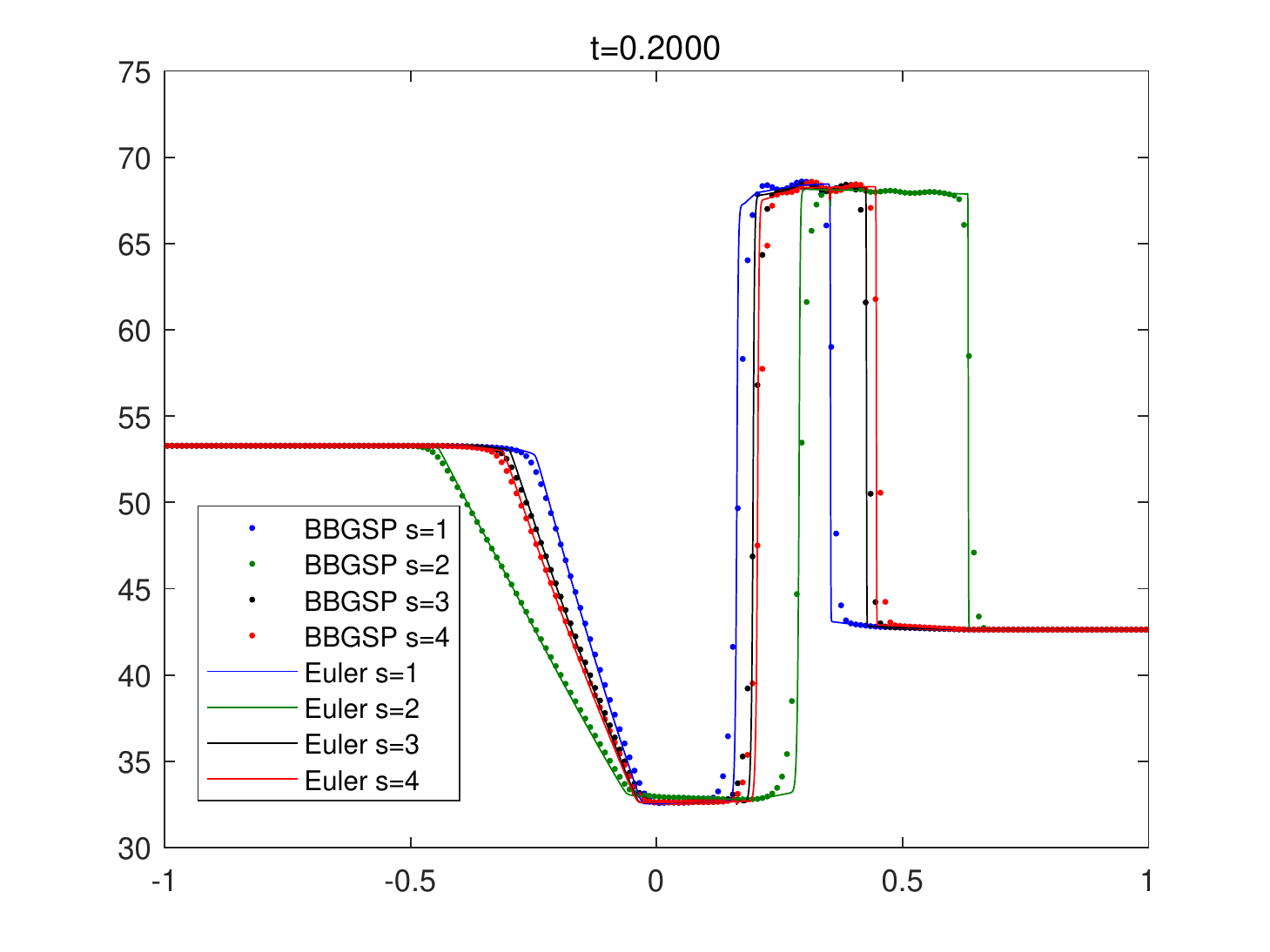}
		\subcaption{Temperature $T_s$, $s=1,2,3,4$}
	\end{subfigure}
	\begin{subfigure}[b]{0.43\linewidth}
		\includegraphics[width=1\linewidth]{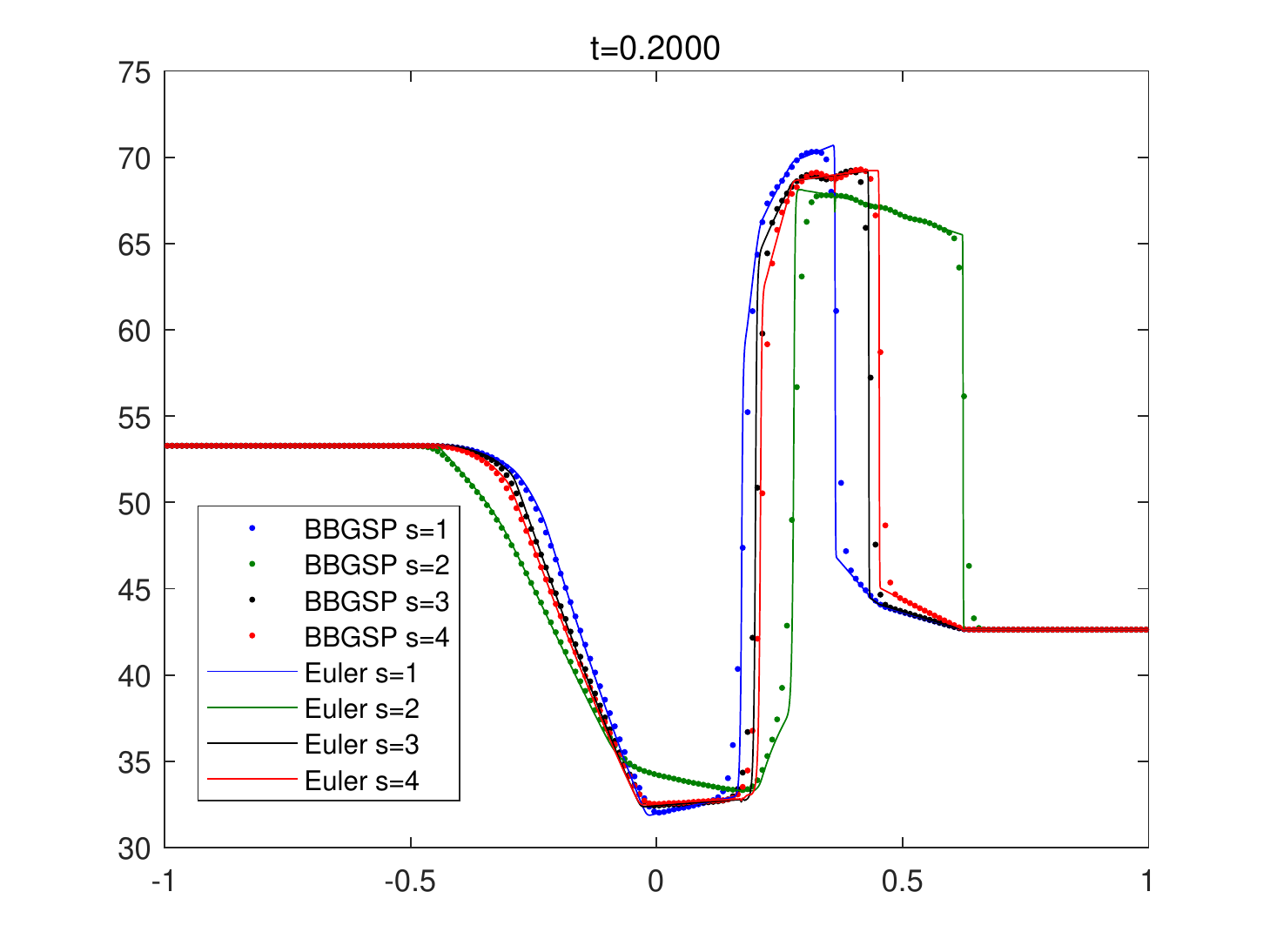}
		\subcaption{Temperature $T_s$, $s=1,2,3,4$}
	\end{subfigure}
	\caption{BDF3-QCW35 for  $(\varepsilon,\kappa)=(10^{-6},10^{-1})$ (left) and $(\varepsilon,\kappa)=(10^{-6},10^{-2})$ (right) associated to initial data in section \ref{sec euler multi}. 
	}\label{fig euler multi 1 0}
\end{figure}

\begin{figure}[htbp]
	\centering
	\begin{subfigure}[b]{0.43\linewidth}
		\includegraphics[width=1\linewidth]{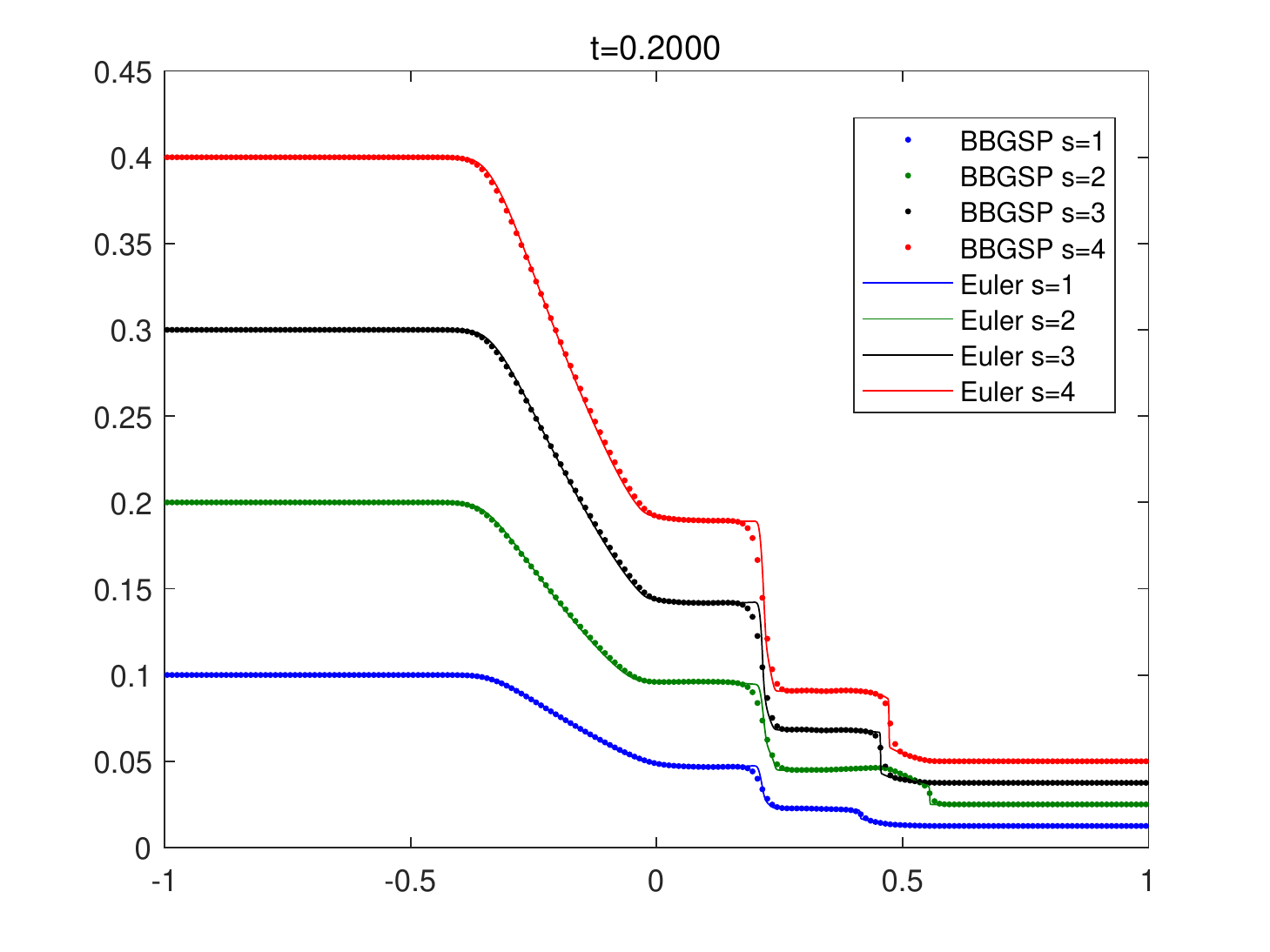}
		\subcaption{Density $\rho_s$, $s=1,2,3,4$}
	\end{subfigure}
	\begin{subfigure}[b]{0.43\linewidth}
		\includegraphics[width=1\linewidth]{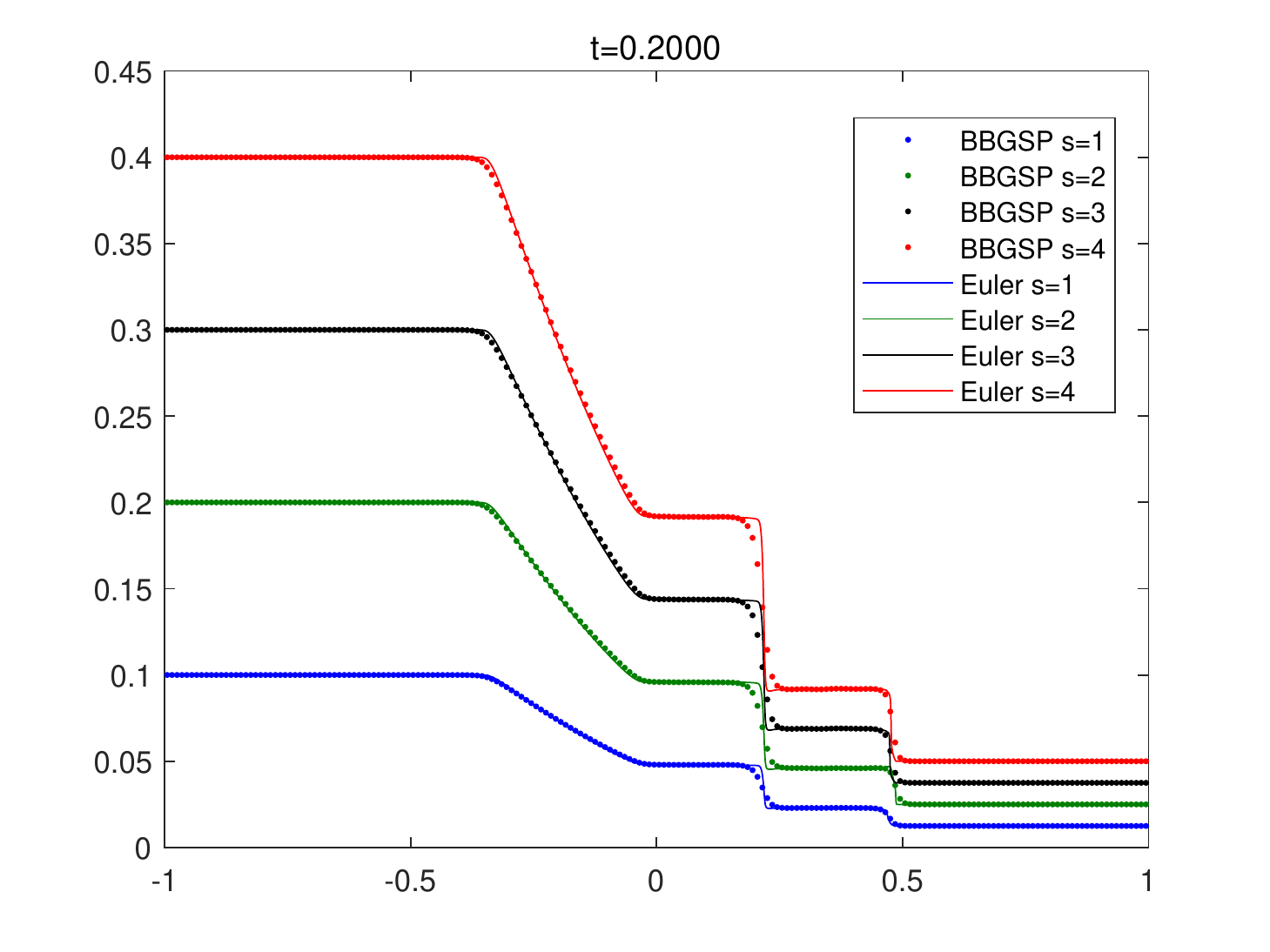}
		\subcaption{Density $\rho_s$, $s=1,2,3,4$}
	\end{subfigure}
\begin{subfigure}[b]{0.43\linewidth}
	\includegraphics[width=1\linewidth]{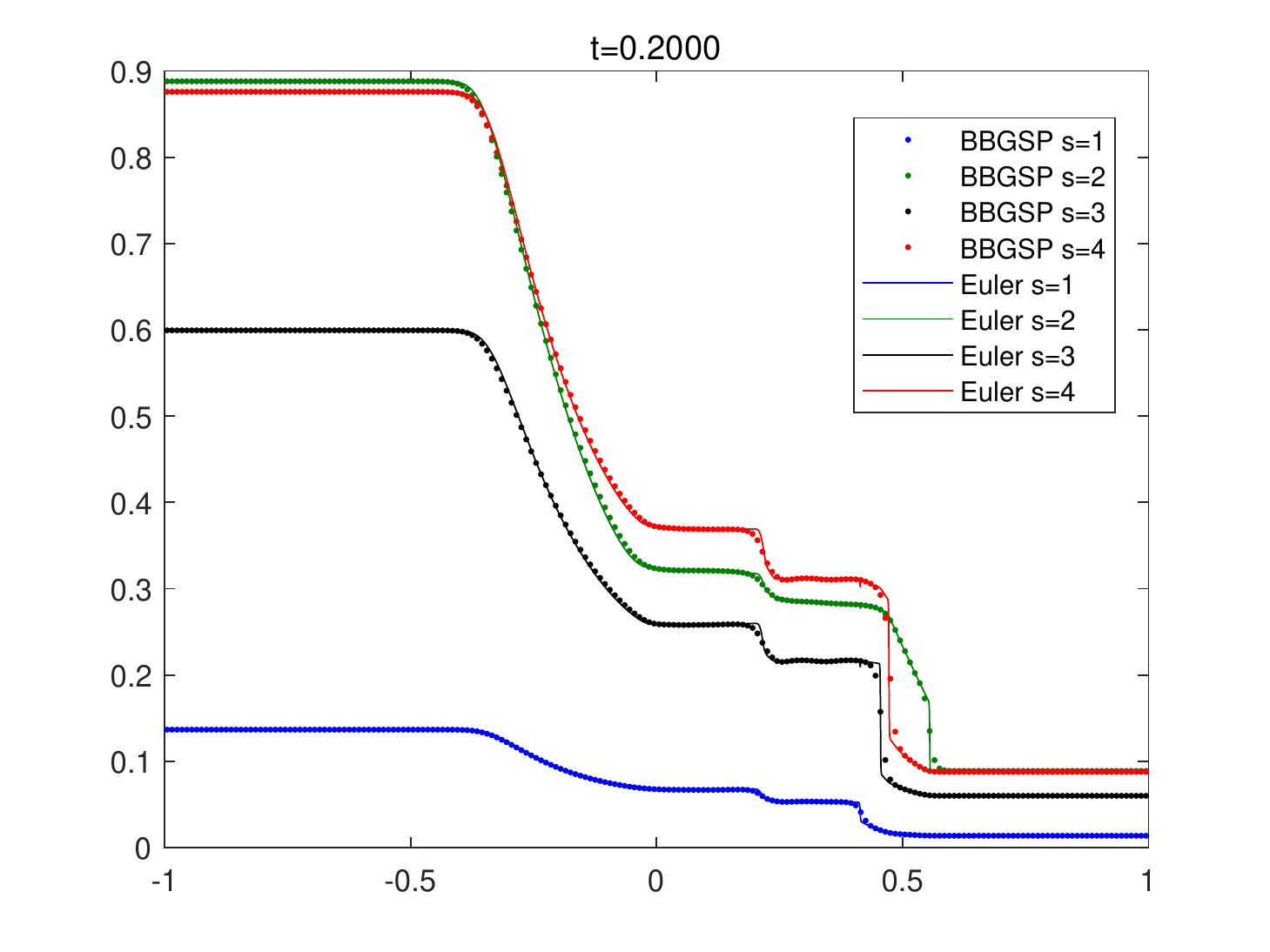}
	\subcaption{Energy $E_s$, $s=1,2,3,4$}
\end{subfigure}
\begin{subfigure}[b]{0.43\linewidth}
	\includegraphics[width=1\linewidth]{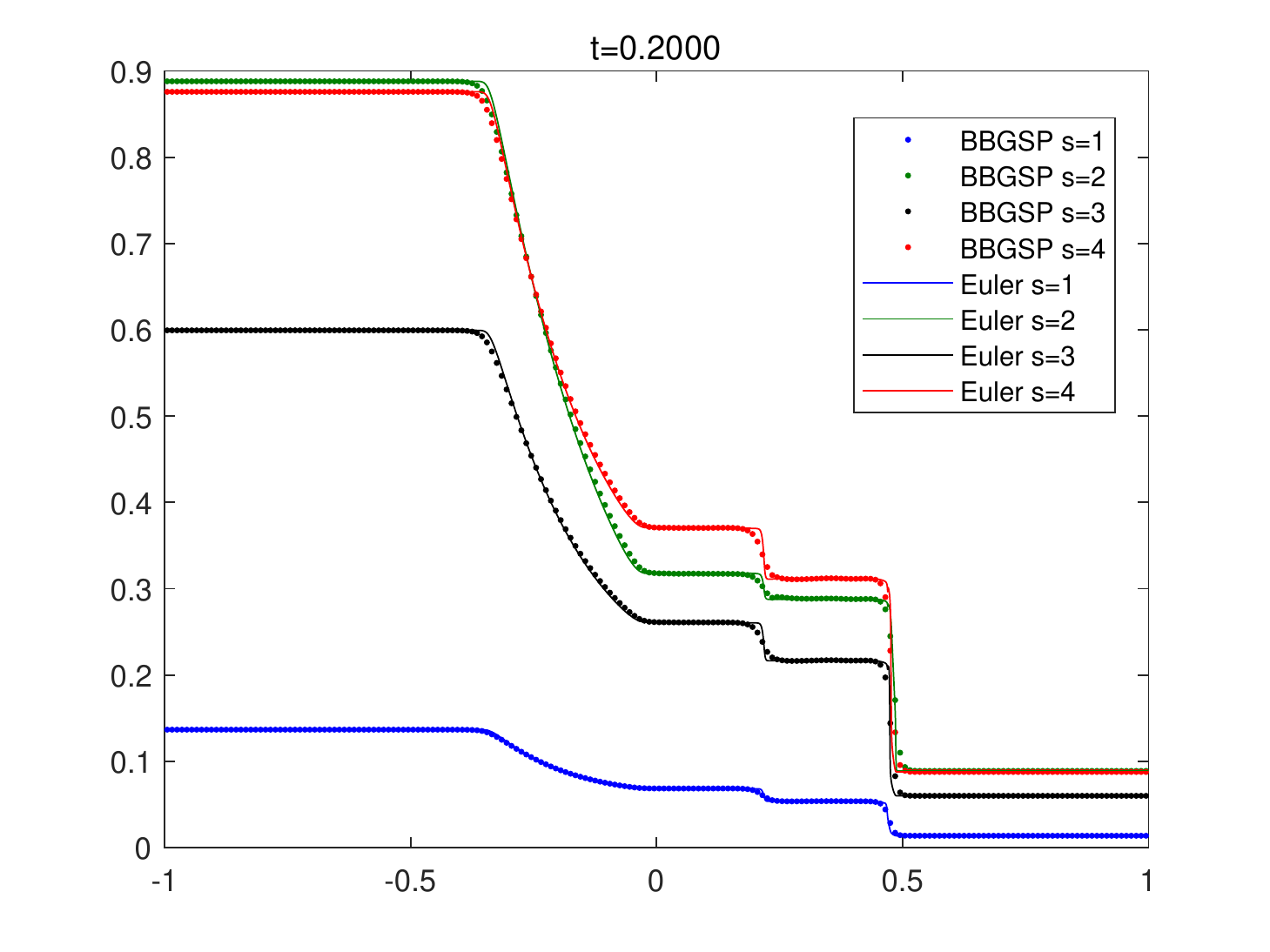}
	\subcaption{Energy $E_s$, $s=1,2,3,4$}
\end{subfigure}
	\begin{subfigure}[b]{0.43\linewidth}
		\includegraphics[width=1\linewidth]{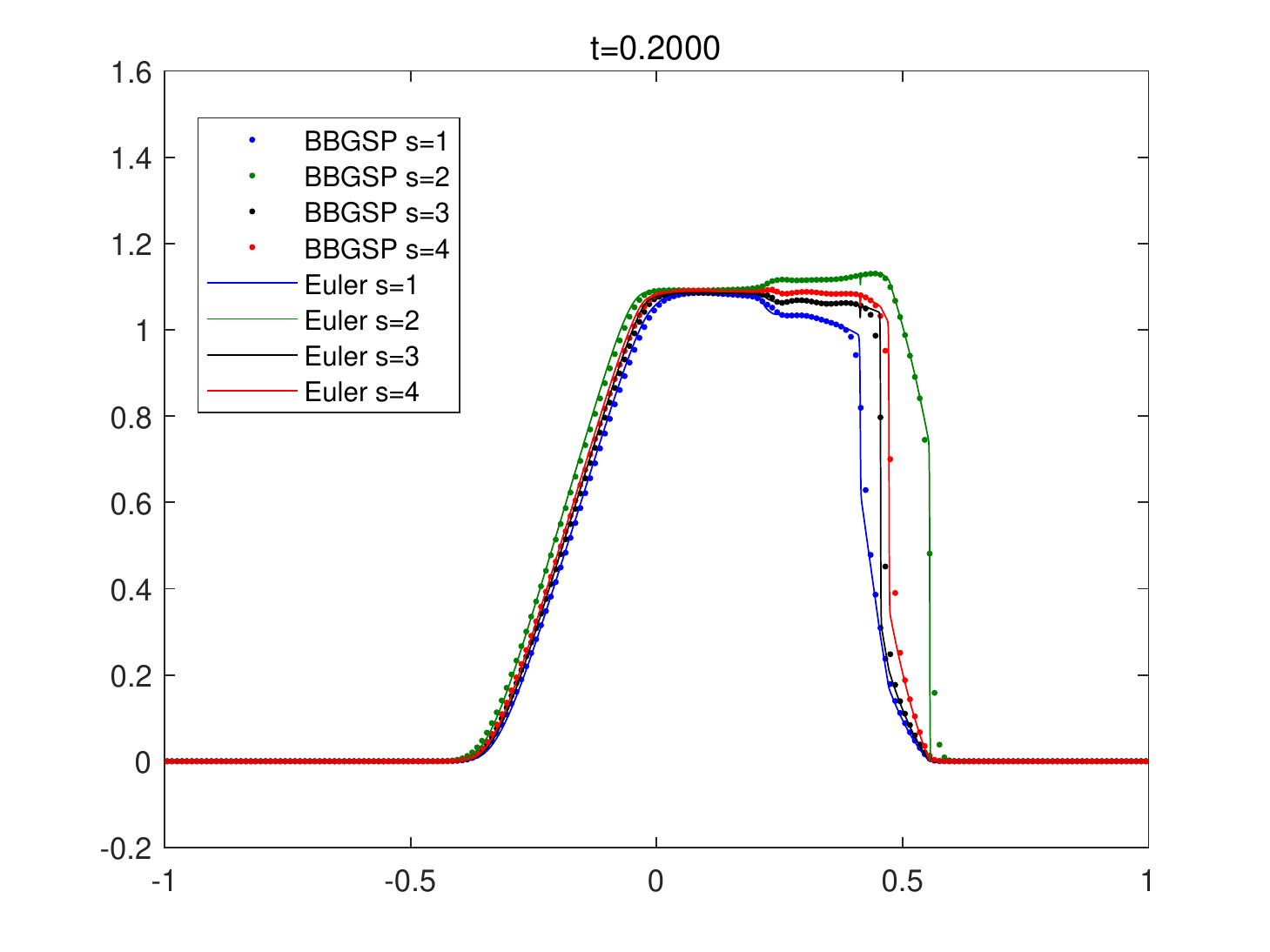}
		\subcaption{Velocity $u_s$, $s=1,2,3,4$}
	\end{subfigure}
	\begin{subfigure}[b]{0.43\linewidth}
		\includegraphics[width=1\linewidth]{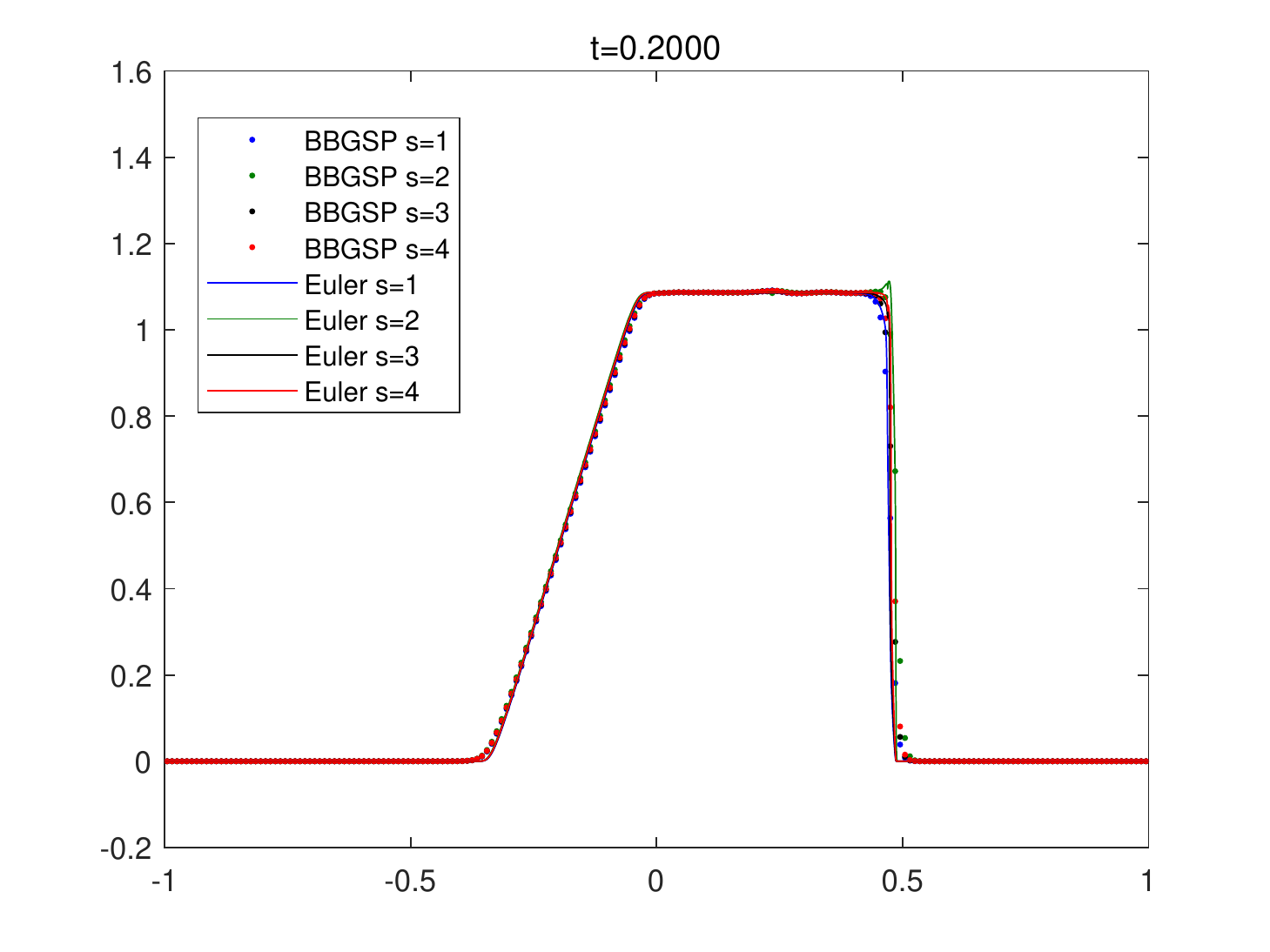}
		\subcaption{Velocity $u_s$, $s=1,2,3,4$}
	\end{subfigure}
	\begin{subfigure}[b]{0.43\linewidth}
		\includegraphics[width=1\linewidth]{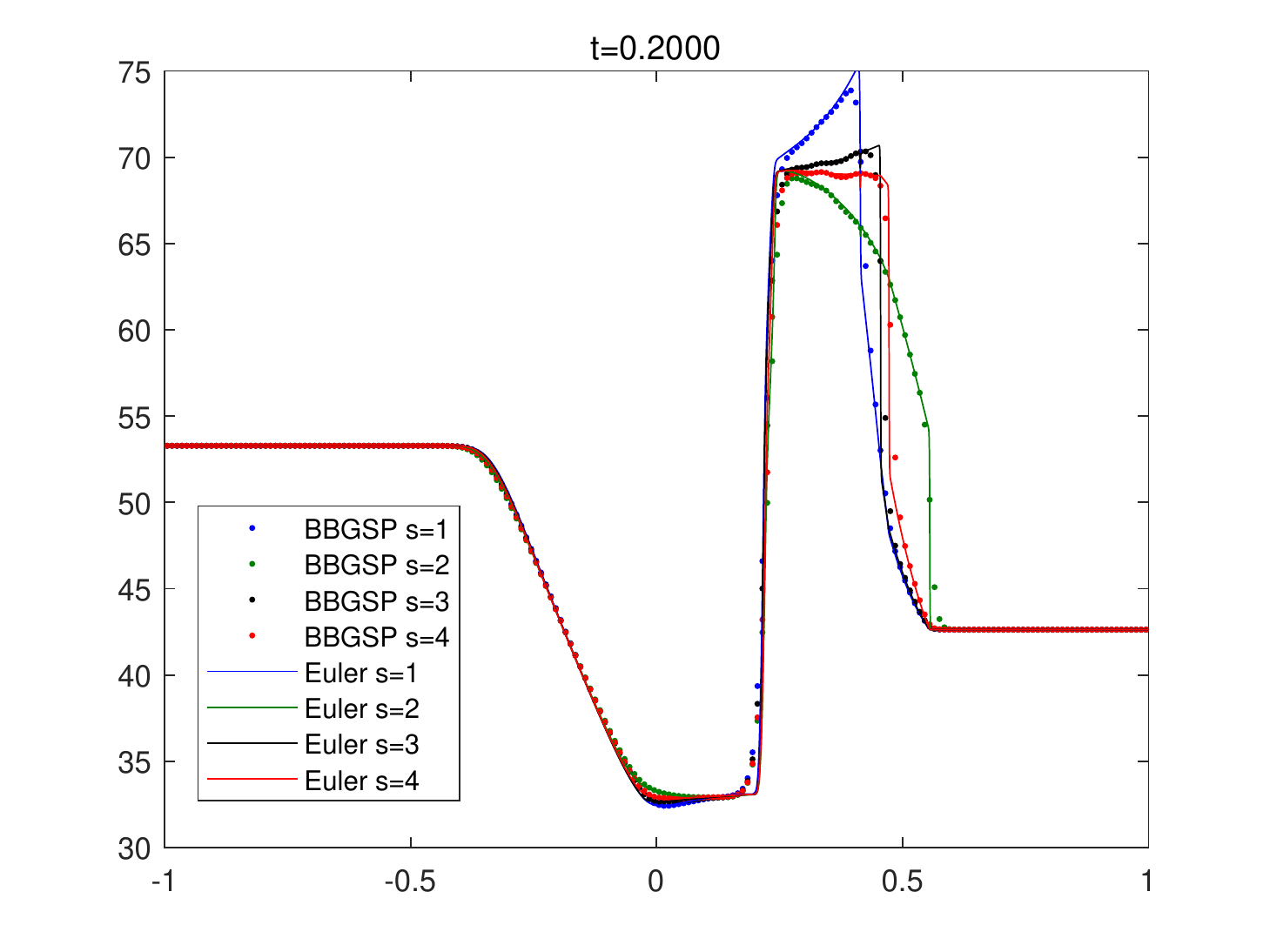}
		\subcaption{Temperature $T_s$, $s=1,2,3,4$}
	\end{subfigure}
	\begin{subfigure}[b]{0.43\linewidth}
		\includegraphics[width=1\linewidth]{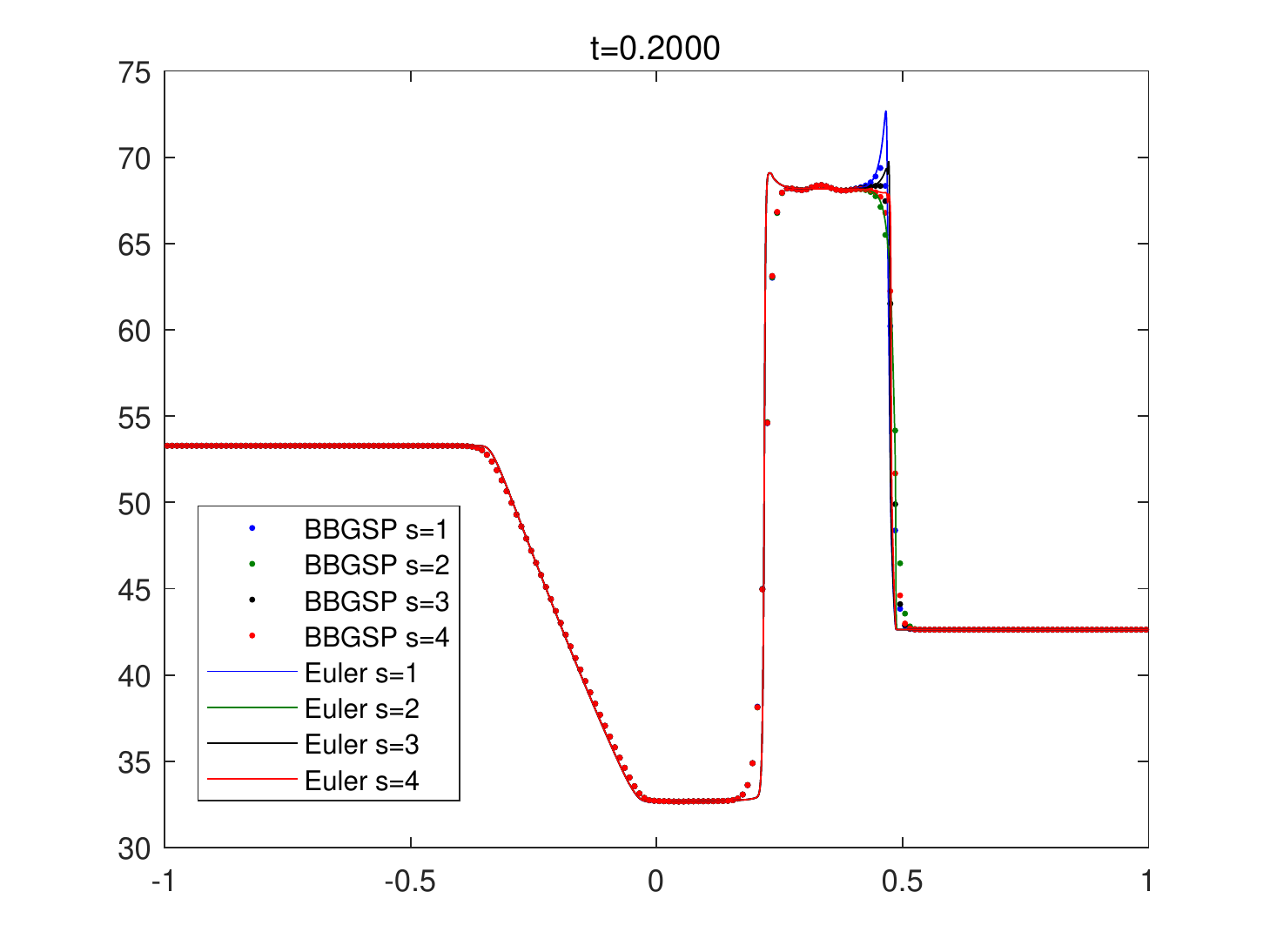}
		\subcaption{Temperature $T_s$, $s=1,2,3,4$}
	\end{subfigure}
	\caption{BDF3-QCW35 for  $(\varepsilon,\kappa)=(10^{-6},10^{-3})$ (left) and $(\varepsilon,\kappa)=(10^{-6},10^{-4})$ (right) associated to initial data in section \ref{sec euler multi}. 
	}\label{fig euler multi -1-2}
\end{figure}

In Figures \ref{fig euler multi 1 0}-\ref{fig euler multi -1-2}, we compare our solutions (dotted lines) with the reference solution (solid lines) to \eqref{Euler multi 1d} obtained by explicit RK3 method with WENO23 with 4000 grid points.

	We first consider two extreme cases: $\kappa=10^{-1}$ (left column of Fig.~\ref{fig euler multi 1 0}) and $\kappa=10^{-6}$ (Fig.~\ref{fig euler}). For $\kappa=10^{-1}$, there are few collisions between gas molecules of different species, and the behavior of each gas is almost independent from the other gases. On the contrary, when there are more frequent collisions $\kappa=10^{-6}$ (right column of Fig.~\ref{fig euler}) each species velocity and temperature converges to the same one. This is what we expect to show with two different Euler models \eqref{Euler 1d} and \eqref{Euler multi 1d}. 
Also, in the intermediate regime, say $\kappa=10^{-3}$, (left column of Fig.~\ref{fig euler multi -1-2}) there is an interesting observation about the species velocity and temperature profiles. 
Gas particles travel from right to left. Light particles are affected by the shock earlier than heavy particles, and in the interaction they tend to change their velocity more quickly. This explains why the change in the energy appears first in the light particles than in the heavy particles. Also, the light particles change their velocities much quicker than the heavy particles, and this is why there is a peak in the kinetic energy, illustrated by the overshoot in the green velocity profile; on the contrary, the effect of the collision on the heavy particles is slightly delayed, and it first manifests with the increase of their temperature, similar to the Brownian motion effect, which explains the blue peak in the temperature profile. All such effects appear on a much shorter time scale for a smaller value of the parameter $\kappa$, as illustrated in Fig.~\ref{fig euler multi -1-2}. From the matches of the numerical solutions and reference solutions for various range of $\kappa$,
	we verified that our SL scheme is able to capture the correct limit to the multi velocities and temperatures Euler system \eqref{Euler multi 1d}.
	We can also confirm that the behavior of the single velocity and temperature Euler system \eqref{Euler 1d}
	can be described by the multi velocities and temperatures Euler system \eqref{Euler multi 1d} in the limit $\kappa \rightarrow 0$.
	Consequently, we can say that our scheme is asymptotic preserving.

\section*{Conclusion}
This papers presents high order conservative semi-Lagrangian schemes applied to a recent consistent BGK model for inert mixtures \cite{BBGSP1}. Thanks to its structure, this BGK model allows different asymptotic limits in collision dominated regimes, depending on if all collisions are equally dominant or if only the intra-species ones are the leading process in the evolution of the mixture. Correspondingly, at zero order asymptotic with respect to a small Knudsen number $\epsilon$, Euler equations involving global velocity and temperature or species velocities and temperatures can be derived. High order conservative semi-Lagrangian schemes presented here are unconditionally stable, and hence it enables to reproduce the correct behavior of Euler equations without resorting to CFL-type time step restriction. The accuracy of the methods and their convergence rates are tested on a sample case. Then, the indifferentiability principle, which can be proved analytically, is shown to be fulfilled  numerically at the desired accuracy. Finally, the AP property is numerically demonstrate on a Riemann problem considering first the classical Euler equations, when all collisions are dominant, and then versus the multi-velocity and temperature Euler system. In this last case, we have shown the trends towards the classical limit, when the parameter $\kappa$ becomes of the same order of $epsilon$. It has been observed the  opposite behaviours in the variation of species velocities and temperature of light and heavy gases, in case of disparate masses. Realistic applications to noble gases, as well as comparisons at Navier Stokes levels, will be discussed in a forthcoming paper \cite{CBGR}.


\section*{Acknowledgement}

S. Y. Cho has been supported by ITN-ETN Horizon 2020 Project ModCompShock, Modeling and Computation on Shocks and Interfaces, Project Reference 642768. S. Y. Cho, S. Boscarino and G. Russo would like to thank the Italian Ministry of Instruction, University and Research (MIUR) to support this research with funds coming from PRIN Project 2017 (No. 2017KKJP4X entitled ``Innovative numerical methods for evolutionary partial differential equations and applications"). 
S. Boscarino has been supported by the University of Catania (``Piano della Ricerca 2016/2018, Linea di intervento 2"). S. Boscarino and G. Russo are members of the INdAM Research group GNCS. M. Groppi thanks the support by the
University of Parma, by the Italian National Group of Mathematical
Physics (GNFM-INdAM), and by the Italian National Research Project ``Multiscale phenomena in Continuum Mechanics:
singular limits, off-equilibrium and transitions" (PRIN 2017YBKNCE).



\end{document}